\documentclass[a4paper]{amsproc}
\usepackage{amsfonts}
\usepackage{amssymb}
\usepackage{amscd}
\usepackage{verbatim}
\usepackage{graphicx}
\usepackage{subfigure}
\usepackage{fancyhdr}
\usepackage[T1]{fontenc}
\usepackage[latin1]{inputenc}
\usepackage[english]{babel}

\input xy
\xyoption{all}

\newenvironment{demo}{\noindent {\it\bf Proof.} \rm}

\newcommand{\fait}[3]{\begin{#1}\label{#2}{#3}\end{#1}}
\newcommand{\cqfd}{\hfill \mbox{$\square$}}

\def\Hom{\mathop{\rm Hom}\nolimits}
\def\rad{\mathop{\rm rad}\nolimits}
\def\End{\mathop{\rm End}\nolimits}
\def\ann{\mathop{\rm ann}\nolimits}

\def\pd{\mathop{\rm pd}\nolimits}
\def\id{\mathop{\rm id}\nolimits}

\def\path{\mathop{\rightsquigarrow}\nolimits}
\def\mod{\mathop{\rm mod}\nolimits}

\def\ind{\mathop{\rm ind}\nolimits}
\def\add{\mathop{\rm add}\nolimits}

\def\pres{\mathop{\rm pres}\nolimits}
\def\fact{\mathop{\rm \in\!\shortmid}\nolimits}

\newcommand{\R}{A[G]}
\newcommand{\F}[1]{F(#1)}
\newcommand{\Ho}[1]{H(#1)}
\newcommand{\Ga}{\Gamma}
\newcommand{\GA}{\Ga(\mod A)}
\newcommand{\mc}{\mathcal}

\def\la{{\mathcal L}_A}
\def\ra{{\mathcal R}_A}
\def\lb{{\mathcal L}_B}
\def\rb{{\mathcal R}_B}
\def\lr{{\mathcal L}_{\R}}
\def\rr{{\mathcal R}_{\R}}

\newcommand{\mor}[3]{\xymatrix@1@C=15pt{#1: #2\ar[r]& #3}}

\begin{document}

\title[Almost laura algebras]{Almost laura algebras \footnote{2000 Mathematics Subject Classification:
16G70, 16G30, 18G05, 16E10.} \footnote{Key words: almost laura
algebras, infinite radical, laura algebras, left and right
supported algebras, skew group algebras}}
\date{}
\author{David Smith}%

\begin{abstract}
\noindent In this paper, we propose a generalization for the class
of laura algebras of \cite{AC03} and \cite{RS04}, which we call
almost laura.  We show that this new class of algebras retains
most of the essential features of laura algebras, especially
concerning the important role played by the non-semiregular
components in their Auslander-Reiten quivers.  Also, we study more
intensively the left supported almost laura algebras, showing that
these are characterized by the presence of a generalized standard,
convex and faithful component. Finally, we prove that almost laura
algebras behave well  with respect to full subcategories,
split-by-nilpotent extensions and skew group algebras.
\end{abstract}

\maketitle

\bibliographystyle{plain}

%
%


In the representation theory of algebras, a prevalent technique
consists of modifying certain features of a well-known family of
algebras in order to obtain one whose representation theory is, to
a large extent, predictable. For instance, in \cite{HRS96},
Happel, Reiten and Smal{\o} defined the quasitilted algebras (that
is the endomorphism algebras of tilting objects over a hereditary
abelian category), thus obtaining a common treatment of both the
classes of tilted and canonical algebras. To overcome some
difficulties caused by the categorical language, they introduced
the left and the right parts of the module category of an algebra
$A$, respectively denoted $\la$ and $\ra$.
%
%
They showed that an algebra $A$ is quasitilted if and only if its
global dimension is at most two and any indecomposable $A$-module
lies in $\la\cup\ra$.

Since then, many generalizations of quasitilted algebras, based on
the behavior of $\la$ and $\ra$ have appeared, such as the shod,
the weakly shod, the laura and the supported algebras (see the
survey \cite{ACLST05}). Among them, laura algebras have been
introduced independently by Assem and Coelho \cite{AC03} and
Reiten and Skowro\'nski \cite{RS04} as a generalization of
representation-finite algebras and weakly shod algebras. Their
nice properties have made them rather interesting and widely
investigated (see \cite{Sk03,AC04,ALR07,LS06,Smith04,DS06}, for
instance).
The aim of this paper is to introduce a new class of algebras,
called almost laura, determined by the behavior of the infinite
radical of $\mod A$ and generalizing laura algebras.

This paper is organized as follows.  In Section
\ref{Preliminaries}, we fix the terminology and prove some
preliminary results.  In Section \ref{Almost laura algebras :
definition and examples}, we give the definition of almost laura
algebras and discuss examples. In Section \ref{Those almost laura
algebras which are laura}, we study the Auslander-Reiten quiver of
an almost laura algebra and we classify the almost laura algebras
which are laura. Section \ref{Supported almost laura algebras} is
devoted to the left (or right) supported almost laura algebras (in
the sense of \cite{ACT04}). Our main result (see (\ref{thm almost
laura})) is an analogue of the result of \cite[(3.1)]{RS04} for
laura algebras (see also \cite[(4.2.5)]{LS06}), and states that if
$A$ is left (or right) supported, then $A$ is almost laura if and
only if its Auslander-Reiten quiver has a generalized standard,
convex and faithful component. Finally, in Section \ref{Full
subcategories, split-by-nilpotent extensions and skew group
algebras}, we show that almost laura algebras behave well with
respect to some constructions preserving homological properties,
such as dealing with full subcategories, split-by-nilpotent
extensions and skew group algebras. The main result of this
section states that if $G$ is a finite group acting on an algebra
$A$ and whose order is invertible in $A$, then $A$ is almost laura
if and only if so is the skew group algebra $A[G]$ (see (\ref{thm
skew})). As a consequence, we get that the infinite radical of $A$
is nilpotent if and only if so is the infinite radical of $A[G]$,
and in this case, they have the same index of nilpotency. We also
deduce that $A$ is cycle-finite (in the sense of \cite{AS90}) if
and only if so is $A[G]$ (see (\ref{cycle-finite})).

%
%
\section{Preliminaries}
    \label{Preliminaries}

In this paper, all algebras are artin algebras over an artinian
ring $k$ (and, unless otherwise specified, connected and basic).
For an algebra $A$, we denote by mod$A$ its category of finitely
generated left modules and by ind$A$ a full subcategory of $\mod
A$ consisting of one representative from each isomorphism class of
indecomposable modules. For a subcategory $\mc{C}$ of $\mod A$, we
write $M\in\mc{C}$ to express that $M$ is an object in $\mc{C}$,
and denote by $\add \mc{C}$ the full subcategory of $\mod A$
having as objects the direct sums of indecomposable summands of
objects in $\mc{C}$. For an $A$-module $M$, we denote by $\pd M$
its projective dimension and by $\id M$ its injective dimension.

We denote by $\GA$ the Auslander-Reiten quiver (AR-quiver for
short) of $A$ and by $\tau_A$
the usual AR-translation. By an \textbf{AR-component} $\Ga$ of
$\Gamma($mod$A)$, we mean a connected component of
$\Gamma($mod$A)$. Then $\Ga$ is \textbf{non-semiregular} if it
contains a projective module and an injective module, and
\textbf{semiregular} otherwise. Also, $\Ga$ is \textbf{faithful}
if it contains a \textbf{faithful} module, that is a module $M$
which cogenerates $A$. Finally, an indecomposable module $M\in
\Ga$ is \textbf{left stable} if $\tau^nM\neq 0$ for each $n\geq 0$
and we define the \textbf{left stable part} of $\Ga$ to be the
full subquiver of $\Ga$ consisting of the left stable modules in
$\Ga$. We define dually the \textbf{right stable} modules and the
\textbf{right stable part} of $\Ga$.

We call \textbf{radical} of mod$A$ and we denote by rad(mod$A$)
the ideal in mod$A$ generated by all non-isomorphisms between
indecomposable modules.  The \textbf{infinite radical}
rad$^\infty(\textrm{mod}A)$ of mod$A$ is the intersection of all
powers rad$^n(\textrm{mod}A)$, with $n\geq 1$, of
rad$(\textrm{mod}A)$. A component $\Ga$ of $\GA$ is
\textbf{generalized standard} \cite{Skow94II} if
$\rad^\infty(M,N)=0$ for each $M,N\in \Ga$.

A \textbf{path of length} $t$ is a sequence %
$ \delta : \xymatrix@C=15pt{M=M_0 \ar[r]^-{f_1} & M_1
\ar[r]^-{f_2} & \cdots \ar[r]^-{f_t} & M_t = N}$
$(t\geq 0)$ where $M_i \in \mbox{ind}A$ and $f_i$ is a non-zero
morphism for each $i$. In this case, we write
$\xymatrix@1@C=15pt{M \ar@{~>}[r]^\delta & N}$
and we say that $M$ is a \textbf{predecessor} of $N$ and $N$ is a
\textbf{successor} of $M$. Following \cite{Skow94II}, a path
$\delta$ is \textbf{infinite} if $f_i\in\rad^\infty(\mod A)$ for
some $i$, and \textbf{finite} otherwise. If each $f_i$ is
irreducible, $\delta$ is a \textbf{path of irreducible morphisms}
and, in this case, $\delta$ is \textbf{sectional} if it contains
no triple $(M_{i-1}, M_i, M_{i+1})$ such that $\tau_A M_{i+1} =
M_{i-1}$. A \textbf{refinement} of $\delta$ is a path
$\xymatrix@C=15pt{M=M'_0 \ar[r]^-{f'_1} & M'_1 \ar[r]^-{f'_2} &
\cdots \ar[r]^-{f'_s} & M'_s = N,}$
with $s\geq t$, with an injective order-preserving function
$\mor{\sigma}{\{1,\dots,t-1\}}{\{1,\dots, s-1\}}$ such that
$M_i=M'_{\sigma(i)}$ when $1\leq i\leq t-1$. Finally, a path
$\delta$ is a \textbf{cycle} if $M=N$ and at least one $f_i$ is
not an isomorphism. An $A$-module $M$ is \textbf{directing} if it
does not lie on any cycle and a component $\Gamma$ of $\GA$ is
\textbf{directed} if it contains only directing modules. Also,
$\Ga$ is \textbf{almost directed} if it contains only finitely
many non-directing modules, and \textbf{quasi-directed} if it is
also generalized standard. Moreover, $\Gamma$ is \textbf{convex}
if any path from $M$ to $N$, with $M,N$ in $\Gamma$, contains only
modules from $\Gamma$.


Let $A$ be an artin algebra.  Following \cite{HRS96}, we define
the left part $\la$ and the right part $\ra$ of $\mod A$ as
follows:
    $$\begin{array}{c}
      \la=\{\ M\in \mbox{ind}A \ | \ \mbox{pd}_A N\leq 1 \mbox{ for each
predecessor } N \mbox{ of } M\ \}, \\
      \ra=\{\ M\in \mbox{ind}A \ | \ \mbox{id}_A N\leq 1 \mbox{ for each
successor } N \mbox{ of } M\ \}. \
    \end{array}$$

The next result is helpful to detect the modules which lie in
$\la$ or in $\ra$.

\fait{lemma}{AC03 (1.6)}{\emph{\cite[(1.6)]{AC03}} \ Let $A$ be an algebra.%
\begin{enumerate}
\item[\emph{(a)}] $\la$
consists of the modules $M\in\ind A$ such that, if there exists a
path from an indecomposable injective module to $M$, then this
path can be refined to a path of irreducible morphisms, and any
such refinement is sectional.
\item[\emph{(b)}]  $\ra$ consists of the modules $N\in\ind A$ such
that, if there exists a path from $N$ to an indecomposable
projective module, then this path can be refined to a path of
irreducible morphisms, and any such refinement is sectional. \cqfd
\end{enumerate}}

We conclude this section with some preliminary results, needed
later on.

\fait{lemma}{pred proj}{Let $A$ be an algebra and $\Ga$ be a
component of $\Ga(\emph{mod} A)$. Assume that
$\rad^\infty(M,N)\neq 0$ for some indecomposable modules $M,N$
with $N\in \Ga$.  Then, for each $L\in \Ga$,
there exists $N'\in\Ga$ such that :%
\begin{enumerate}
\item[\emph{(a)}] There exists a path of irreducible morphisms
from $N'$ to $N$;
\item[\emph{(b)}] $N'$ is a predecessor of $L$ or is a predecessor
of a projective module in $\Ga$;
\item[\emph{(c)}] $\rad^\infty(M,N')\neq 0$.
\end{enumerate}}
\begin{demo}
Let $M$ and $N$ be as in the statement. There exists a path of
infinite length of irreducible morphisms
$$\xymatrix@1@C=20pt{\cdots \ar[r] & N_r\ar[r]^{h_r} &
N_{r-1}\ar[r] & \cdots \ar[r]^{h_2} & N_1 \ar[r]^-{h_1} & N_0=N}$$
in $\ind A$ such that there exists $u_r\in \rad^\infty(M,N_r)$
with $h_1h_2\cdots h_ru_r\neq 0$ for each $r\geq 1$ (see
\cite[(2.1)]{Sk94II}). Let $L\in\Ga$. We claim that there exists
$s\geq 1$ such that $N_s$ is a predecessor of $L$ or is a
predecessor of a projective module in $\Ga$. Indeed, if this is
not the case, then $N_i$ is not projective for all $i$ and it
follows from \cite[(1.1)]{CL02} that there exists an integer
$r\geq 1$ which is minimal for the property that $N_i$ is not a
predecessor of $\tau^rN$ for all $i$. By the choice of $r$, there
exists $N_j$ such that $N_j$ is a predecessor of $\tau^{r-1}N$. We
claim that the path
$\xymatrix@1@C=15pt{N_m\ar[r]^{h_m} & N_{m-1}\ar[r]^{h_{m-1}} &
\cdots \ar[r]^{h_{j+1}} & N_j}$
is sectional for each $m>j$. Indeed, if this is not the case, then
there exists $n$ with $j\leq n \leq m-2$ such that $N_{n+2}=\tau
N_n$. This yields a path
$\xymatrix@1@C=15pt{N_{n+2}=\tau N_n\ar@{~>}[r] & \tau N_j
\ar@{~>}[r] & \tau^rN},$
a contradiction to the choice of $r$. In particular, $N_m\neq N_n$
whenever $m\neq n$ and $m,n\geq j$. Therefore, $\Hom(N_m, \tau
N_n)\neq 0$ for some $m,n\geq j$ by \cite[(Lemma 2)]{Sk94III}.
Again, this yields a path from $N_m$ to $\tau^rN$,
%
%
a contradiction. Thus there exists $s\geq 1$ such that $N_s$ is a
predecessor of $L$ or is a predecessor of a projective in $\Ga$.
\cqfd
\end{demo}

As immediate consequences, we obtain the following corollary which
generalizes results obtained in \cite[(1.4)]{AC03} and
\cite[(1.4)]{Smith04}.

\fait{cor}{paths to cycles}{Let $A$ be an algebra, $\Ga$ be a
component of $\Ga(\emph{mod}A)$ and assume that $M$ is a
non-directing module
in $\Ga$.%
\begin{enumerate}
\item[\emph{(a)}] If $\Ga$ contains projective modules, then there exists
a path from $M$ to a projective module in $\Ga$.
\item[\emph{(b)}] If $\Ga$ contains injective modules, then there exists
a path from an injective module in $\Ga$ to $M$.
\end{enumerate}}
\begin{demo}
We only prove (a) since the proof of (b) is dual.\\ %
(a). Let
$\xymatrix@1@C=15pt{M=M_0 \ar[r]^-{f_1} & M_1\ar[r]^{f_2} & \cdots
\ar[r]^-{f_t} & M_t=M}$
be a cycle in $\ind A$.  If no $f_i$ belongs to $\rad^\infty(\mod
A)$, then this cycle can be refined to a cycle of irreducible
morphisms in $\Ga$, and the result follows from
\cite[(1.4)]{AC03}.  Otherwise, we have
$f_{i}\in\rad^\infty(M_{i-1}, M_i)$ for some $M_i\in \Ga$, and it
follows from (\ref{pred proj}) that there exists a projective
module $P$ in $\Ga$ and a path from $M_{i-1}$ to $P$. This gives a
path from $M$ to $P$ as required. \cqfd
\end{demo}

We also deduce from (\ref{pred proj}) the following generalization
of \cite[(1.5)]{ACT04}.

\fait{cor}{directed}{Let $A$ be an algebra and $\Ga$ be a
component of $\Ga(\emph{mod}A)$.
\begin{enumerate}
\item[\emph{(a)}] If $\Ga$ contains projectives, then $\ra\cap\Ga$ contains only directing modules.
\item[\emph{(b)}] If $\Ga$ contains injectives, then $\la\cap\Ga$ contains only directing modules.
\end{enumerate}}
\begin{demo}
We only prove (a) since the proof of (b) is dual.\\ %
(a). Assume that $M\in \ra\cap\Ga$ and
$\omega : \xymatrix@1@C=15pt{M \ar@{~>}[r]&M}$
is a cycle in $\ind A$. By (\ref{paths to cycles}), there exists a
path
$\xymatrix@1@C=15pt{M \ar@{~>}[r]^{\omega} & M \ar@{~>}[r] & P}$
where $P$ is projective.  By (\ref{AC03 (1.6)}), this path can be
refined to a sectional path of irreducible morphisms. But this
contradicts the non-sectionality of cycles \cite{BS83,IT84}. \cqfd
\end{demo}



%
%
\section{Almost laura algebras : definition and examples}
    \label{Almost laura algebras : definition and examples}

We recall from \cite{AC03} that an artin algebra $A$ is called
\textbf{laura} if the set $\ind A \setminus (\la\cup\ra)$ is
finite. Since the left and the right part generally behave well,
the spirit of laura algebras is to deal with algebras having
potentially only finitely many "unpredictable" modules.  This idea
behind almost laura algebras is to accept infinitely many such
modules but restrict their scope by adding a condition on the
morphisms between them.

\fait{defin}{Almost laura algebras}{\emph{%
An artin algebra is called \textbf{almost laura} if
$\rad^\infty(M,N)$ vanishes for all $M, N\in \ind A
\setminus(\la\cup\ra)$.}}

In the vein of \cite{ACLST05}, we also say that an almost laura
algebra is \textbf{strict} if it is not quasitilted.
The following proposition provides many equivalent useful
conditions for an algebra to be almost laura.

\fait{prop}{prop almost laura}{Let $A$ be an algebra.  The
following
are equivalent: %
\begin{enumerate}
\item[\emph{(a)}] $A$ is almost laura.
\item[\emph{(b)}] For all $M \in \ind A \setminus \la$ and $N\in \ind A \setminus \ra$, we have
$\rad^\infty(M,N)=0$.
\item[\emph{(c)}] There is no infinite path between modules in $\ind A \setminus(\la\cup\ra)$.
\item[\emph{(d)}] There is no infinite path from a module not in $\la$ to a module not in $\ra$.
\item[\emph{(e)}] There is no infinite path from an injective module to a projective module.
\item[\emph{(f)}] There is no infinite path from a module $M$, with $\pd M\geq 2$, to a module
$N$, with $\id N\geq 2$.
\end{enumerate}}
\begin{demo} The equivalence of (a), (b), (c) and (d) follows
from the fact that $\la$ is closed under predecessors and $\ra$ is
closed under successors.\\
(e) implies (f). Let
$\xymatrix@1@C=15pt{M\ar@{~>}[r]^\omega & N}$
be a path in $\ind A$, with $\pd M \geq 2$ and $\id N\geq 2$.
Since $\pd M \geq 2$, we have $\Hom_A(I, \tau M)\neq 0$ for some
indecomposable injective $I$ and so there exists a path
$\xymatrix@1@C=15pt{\omega' \ : \ I\ar@{~>}[r] & M}$
in $\ind A$. Dually, there exists a path
$\xymatrix@1@C=15pt{\omega'' \ : \ N\ar@{~>}[r] & P}$
for some indecomposable module $P$.  This yields a path
$\xymatrix@1@C=15pt{I\ar@{~>}[r]^-{\omega'} & M
\ar@{~>}[r]^-{\omega} & N \ar@{~>}[r]^-{\omega''} & P}$,
which is finite by assumption, whence so is $\omega$.\\
(f) implies (d). This clearly follows from the definitions of
$\la$ and $\ra$, since any path from a module not in $\la$ to a
module not in $\ra$ can be extended to a path from a module having
projective dimension at least two to a module having injective
dimension at least two. \\
(d) implies (e). Let
$\delta:\xymatrix@1@C=15pt{I=M_0\ar[r]^-{f_1}& M_1\ar[r]^{f_2} &
\cdots \ar[r]^-{f_t} & M_t=P}$ be a path in $\ind A$ from an
injective $I$ to a projective $P$. Assume that $f_i\in
\rad^\infty(\mod A)$, for some $1\leq i\leq t$. For any $n\geq 0$,
it follows from \cite[(1.1)]{Smith04} that $\delta$ may be refined
to a path
$$\delta' :\xymatrix@1@C=18pt{I=M_0\ar@{~>}[r] &
M_{i-1}\ar[r]^{h_0}&N_0\ar[r]^{h_1}&N_1\ar[r]^{h_2} & \cdots
\ar[r]^{h_{n}} & N_{n} \ar[r]^{g_n} & M_{i} \ar@{~>}[r]&P}$$
where $g_n\in \rad^\infty(\mod A)$ and $N_k\neq N_l$ whenever
$k\neq l$. Since there are only finitely many modules in $\la$
which are successors of an injective by \cite[(1.5)]{AC03} (see
also \cite[(3.2.6)]{LS06}), there exists $n\geq 0$ such that
$N_n\notin \la$. Applying the dual argument to $g_n$ yields an
infinite path $\delta'' :\xymatrix@1@C=18pt{N_n \ar@{~>}[r]& M}$,
with $M\notin \ra$, a contradiction to the hypothesis.  \cqfd
\end{demo}

We get the following corollary as an immediate consequence of
(\ref{prop almost laura})(e).

\fait{cor}{cor counter-ex}{If $A$ is an almost laura algebra, then
$\rad^\infty(I, P)=0$ for any injective $A$-module $I$ and
projective $A$-module $P$. \cqfd}

\fait{remark}{counter-ex}{\emph{%
We stress that the converse of the above corollary is false, as
can be easily verified with the radical square zero algebra $A$
given by the quiver
$\xymatrix@1@C=20pt{1 \ar@<.5ex>[r] \ar@<-.5ex>[r] & 2
\ar@<.5ex>[r] \ar@<-.5ex>[r] & 3 \ar@<.5ex>[r] \ar@<-.5ex>[r] &
4}.$}}

We now gives few examples of almost laura algebras.

\fait{examples}{ex almost laura}{ \emph{ %
\begin{enumerate}
\item[(a)] By \cite[(3.3)]{AC03}, any laura algebra is almost laura. In
particular, so is any representation-finite or quasitilted
algebra.
\item[(b)] Let $A$ be the algebra given by the quiver
$$\xymatrix@1@R=2pt{%
1 \ar@<.5ex>[r]^{\beta_1} \ar@<-.5ex>[r]_{\beta_2} & %
2 \ar[dr]^{\alpha}
\\
&& \ 5 \\
3 \ar@<.5ex>[r]^{\delta_1} \ar@<-.5ex>[r]_{\delta_2} & %
4 \ar[ur]_{\gamma}} $$
bound by $\alpha\beta_2 = \gamma\delta_1 = \gamma\delta_2=0$. Then
$\GA$ has the shape presented in Fig.~1 below (where
indecomposable modules are represented by their Loewy series),
where we identify both copies of the module $\substack{2\\ 1}$
along the vertical dashed line, and both copies of the module $2$
along the horizontal dashed line. The horizontal dotted lines
represent the AR-translations. One can verify that $A$ is an
almost laura algebra, but not a laura algebra.
\end{enumerate}}}
\begin{figure}[!htb]
\begin{center}
\begin{picture}(0,0)%
\includegraphics{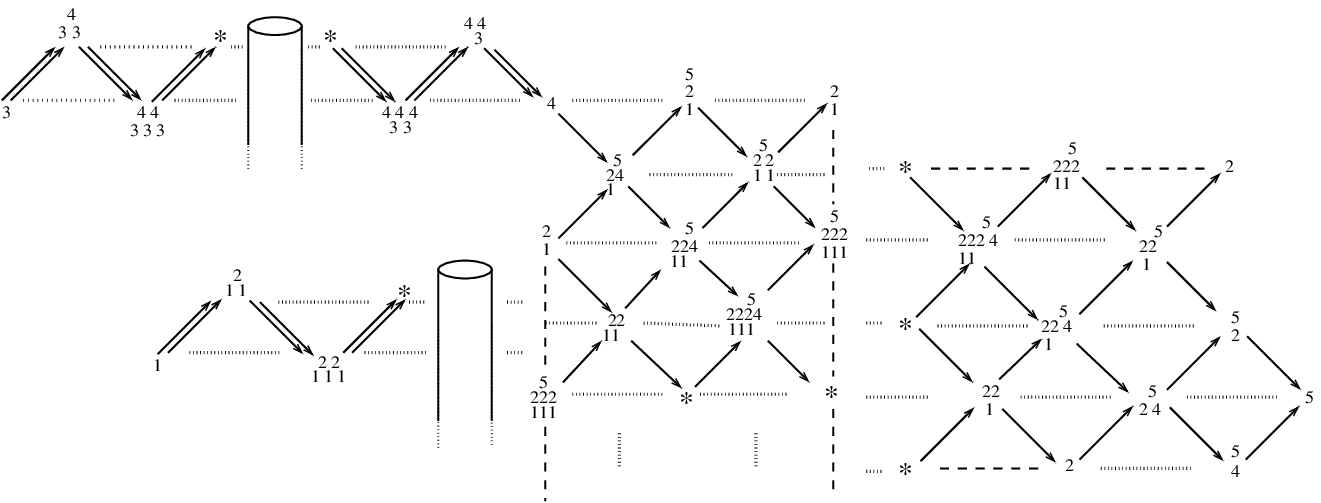}%
\end{picture}%
\setlength{\unitlength}{750sp}%
\begingroup\makeatletter\ifx\SetFigFont\undefined%
\gdef\SetFigFont#1#2#3#4#5{%
  \reset@font\fontsize{#1}{#2pt}%
  \fontfamily{#3}\fontseries{#4}\fontshape{#5}%
  \selectfont}%
\fi\endgroup%
\begin{picture}(33174,12494)(707,-16605)
\end{picture}
    \caption{$\GA$}
    \end{center}
\end{figure}
%

In this latter example, the algebra has been obtained by
performing a one-point extension in a chosen homogeneous tube of
the Kronecher algebra formed by the vertices $1$ and $2$, and by
"gluing" another Kronecker algebra to the resulting ray tube.
Repeating the same procedure in another tube would result in an
almost laura algebra having two non-semiregular components. Since
there are infinitely many such tubes, this shows that one can
construct almost laura algebras having arbitrarily many
non-semiregular components.

We would like to propose the following problem, which is is an
analogue to Skowro{\'n}ski's conjecture for laura algebras
\cite{Sk03}.

\fait{prob}{prob1}{Let $A$ be an algebra. Are the following
conditions equivalent?
\begin{enumerate}
\item[\emph{(a)}] $A$ is almost laura.
\item[\emph{(b)}] $\rad^\infty(M,N)=0$ for all $M,N\in\ind A$, with $\pd M\geq 2$ and
$\id N \geq 2$.
\item[\emph{(c)}] There is no infinite path between modules having both projective and injective dimensions at least $2$.
\end{enumerate}}

%

%
%
\section{Those almost laura algebras which are laura}
    \label{Those almost laura algebras which are laura}

The definition of almost laura algebras is closely related to that
of laura algebras.  In this section, we are interested in
determining when an almost laura algebra is laura. We recall that
strict laura algebras are characterized by the existence of a
unique non-semiregular component in their AR-quiver, which is
moreover quasi-directed and faithful (see \cite{AC03,RS04}).  Our
approach consists in studying the behavior of the non-semiregular
components in the AR-quiver of almost laura algebras.  As we shall
see, those components behave similarly as for laura algebras. We
infer some characterizations of almost laura algebras which are
laura. Our results on the non-semiregular components will also
play a major role in Section \ref{Supported almost laura
algebras}.

\subsection{Non-semiregular components and almost laura algebras}
    \label{Non-semiregular components and almost laura algebras}

We begin our investigation of non-semiregular components over
almost laura algebras with the following key lemma, whose proof is
a routine application of (\ref{prop almost laura}) and (\ref{pred
proj}). We leave the verification to the reader.

\fait{lemma}{lem paths}{An algebra is almost laura if and only if
there is no infinite path from a module $M$ lying in a component
containing injectives to a module $N$ lying in a component
containing projectives.
\cqfd}
%

As a first application, we get the following corollary.

\fait{cor}{cor LAURA}{Let $A$ be an almost laura algebra. If
$M\in\ind A\setminus(\la\cup\ra)$, then $M$ belongs to a
non-semiregular component of $\GA$.}
\begin{demo}
Let $M\notin \la\cup\ra$.  Then there exists a path
$\xymatrix@1@C=15pt{I\ar@{~>}[r]& M\ar@{~>}[r] & P}$
for some injective module $I$ and some projective module $P$.
Since $A$ is almost laura, this path is finite and so $I$ and $P$
belong to the same component as $M$. \cqfd
\end{demo}

We recall from \cite{CS96} that the AR-quiver $\GA$ of a
quasitilted algebra $A$ generally does not contain non-semiregular
components, but if it does, then it contains a unique
non-semiregular component $\Ga$.  Moreover, the algebra $A$ is
then tilted and $\Ga$ is the unique connecting component of $\GA$.
It is well-known that a tilted algebra always admits exactly one
or two connecting components. On the other hand, any strict laura
algebra admits non-semiregular components (see \cite{AC03}). The
following proposition states that the same is true for strict
almost laura algebras.

\fait{prop}{prop nsr}{Let $A$ be an almost laura algebra. %
\begin{enumerate}
\item[\emph{(a)}] If $\Ga$ is a non-semiregular component in $\Ga(\emph{mod} A)$, then $\Ga$ is
generalized standard and convex.%
\item[\emph{(b)}] If $A$ is a strict almost laura algebra, then $\Ga(\emph{mod} A)$ contains non-semiregular components.
\end{enumerate}}
\begin{demo}
(a). This directly follows from the lemma.\\ %
(b). Since $A$ is not quasitilted, it follows from
\cite[(II.1.14)]{HRS96} that there exists an indecomposable
projective module $P$ not lying in $\la$. So, there is a path from
an injective module $I$ to $P$ in $\ind A$.  Since $A$ is almost
laura, the modules $P$ and $I$ belong to the same component of
$\GA$, which is thus non-semiregular. \cqfd
\end{demo}

\fait{remark}{rem repfin}{\emph{%
The above result has a direct nice consequence.  In fact, a
well-known conjecture in representation theory of algebras states
that if an algebra $A$ has a connected AR-quiver, then $A$ is
representation-finite.  Since the AR-quiver then consists of a
unique non-semiregular component, and $A$ is representation-finite
if and only if $\rad^\infty(\mod A)= 0$ by Auslander's theorem
(see \cite[(V.7.7)]{ARS97}), it follows from the above proposition
that the conjecture has a positive answer for almost laura
algebras.  In other words, if $A$ is an almost laura algebra such
that $\GA$ is connected, then $A$ is representation-finite.}}
%

For the remaining part of this section, we let $A$ be an almost
laura algebra and $\Ga$ be a non-semiregular component of $\GA$.
Here and in the sequel, we also use the following notation: if
$\mc{A}$ and $\mc{B}$ are two classes of $A$-modules, then we
write $\Hom_A(\mc{A},\mc{B})\neq 0$ to express that there exists a
non-zero morphism from a module in $\mc{A}$ to a module in
$\mc{B}$.

The following are generalizations of \cite[(4.1)]{AC03} and
\cite[(4.2)]{AC03}. The proof of the lemma follows directly from
(\ref{lem paths}) and it is omitted.

\fait{lemma}{lem I-paths}{Let $A$ and $\Ga$ be as above.
\begin{enumerate}
\item[\emph{(a)}] Assume that $I$ is an indecomposable
injective module such that there exists a path
$\xymatrix@1@C=15pt{I\ar@{~>}[r]& M}$ with $M\in \Ga$, then $I$
belongs to $\Ga$.
\item[\emph{(b)}] Assume that $P$ is an indecomposable
projective module such that there exists a path
$\xymatrix@1@C=15pt{M\ar@{~>}[r]& P}$ with $M\in \Ga$, then $P$
belongs to $\Ga$. \cqfd
\end{enumerate}}
%

\fait{prop}{prop Hom}{Let $A$ and $\Ga$ be as above, and let
$\Ga'$ be a component of $\GA$ distinct from $\Ga$.
\begin{enumerate}
\item[\emph{(a)}] If $\Hom_A(\Ga', \Ga)\neq 0$, then $\Ga' \subseteq\la\setminus\ra$.
\item[\emph{(b)}] If $\Hom_A(\Ga, \Ga')\neq 0$, then $\Ga'\subseteq\ra\setminus\la$.
\item[\emph{(c)}] Either $\Hom_A(\Ga', \Ga)= 0$, or $\Hom_A(\Ga,
\Ga')=0$.
\end{enumerate}}
\begin{demo}
(a). Let $M, M'\in \Ga'$, $N\in\Ga$ and assume that $0\neq f \in
\Hom_A(M,N)$. We need to show that $M'\in\la\setminus\ra$. Clearly
$f\in \rad^\infty(\mod A)$. By (\ref{pred proj}), there exists
$N'\in\Ga$ such that $N'$ is a predecessor of a projective $P$ in
$\Ga$ and $\rad^\infty(M,N')\neq 0$. Dually, there exists $M''\in
\Ga'$ such that $M''$ is a successor of $M'$ or a successor of an
injective module in $\Ga'$ and $\rad^\infty(M'', N')\neq 0$.  By
(\ref{lem paths}), $M''$ is not a successor of an injective.  So
there exists a path
$\xymatrix@1@C=15pt{M'\ar@{~>}[r]& M'' \ar[r]^g & N' \ar@{~>}[r]&
P}$
where $g$ is a non-zero morphism in $\rad^\infty(M'',N')$. Then
$M'\in\la\setminus\ra$ by (\ref{AC03 (1.6)}).  So $\Ga'\subseteq
\la\setminus\ra$.\\
(b). The proof is dual to that of (a).\\%
(c). This follows directly from (a) and (b). \cqfd
\end{demo}

We prove in (\ref{thm Hom2}) below a stronger version of this
result when $A$ is left (or right) supported.
We conclude with an observation on semiregular components.

\fait{prop}{prop semiregular}{Let $A$ be an almost laura algebra
and $\Ga'$ be a semiregular component of $\GA$.
\begin{enumerate}
\item[\emph{(a)}] $\Ga'\subseteq \la\cup\ra$.
\item[\emph{(b)}] If $\Ga'$ contains injectives but no projectives, then
$\Ga'\subseteq \ra$.
\item[\emph{(c)}] If $\Ga'$ contains projectives but no
injectives, then $\Ga'\subseteq \la$.
\item[\emph{(d)}] If $\Ga'$ is regular, that is it contains neither injectives nor projectives, then
$\Ga'$ lies in $\la\setminus\ra$, $\ra\setminus\la$ or
$\la\cap\ra$.
\end{enumerate}}
\begin{demo}
(a). This directly follows from (\ref{cor LAURA}).\\
(b). Assume that $M$ is a module in $\Ga'$ which does belong to
$\ra$. By (\ref{AC03 (1.6)}) there exists a path $\delta$ from $M$
to a projective module $P$. Since $P\notin\Ga'$ by assumption,
this path is infinite.  By the dual of (\ref{pred proj}), there
exists an infinite path from an injective module in $\Ga'$ to $P$,
which contradicts the fact that $A$ is almost laura by (\ref{prop
almost laura}).\\
(c). The proof is dual to that of (b).\\
(d). In view of (a), it suffices to show that if
$\Ga'\cap\la\neq\emptyset$ (or $\Ga'\cap\ra\neq\emptyset$), then
$\Ga'\subseteq\la$ (or $\Ga'\subseteq\ra$ respectively). Assume
that $\Ga'\cap\la\neq\emptyset$ and let $M,N\in\Ga'$ with
$M\in\la$.  If $N\notin\la$, then there exists by (\ref{AC03
(1.6)}) a path $\delta$ from $N$ to an injective module $I$.  But
then, since $I\notin\Ga'$, this path is infinite and it follows
from (\ref{pred proj}) that there exists an infinite path from $M$
to $I$, contradicting the fact that $M\in\la$.  So
$\Ga'\subseteq\la$.  Similarly $\Ga'\cap\ra\neq\emptyset$ implies
$\Ga'\subseteq\ra$. \cqfd
\end{demo}

\subsection{On almost laura algebras which are laura}

In this section, we provide necessary and sufficient conditions
for an almost laura algebra to be laura and also deduce new
characterizations of laura and weakly shod algebras. We begin with
the following key lemma.

\fait{lemma}{lem std gen}{Let $A$ be an algebra and $\Ga$ be a
generalized standard and convex component of $\GA$. For all
$L,N\in \Ga$, there are only finitely many directing modules $M$
lying on a path $\xymatrix@1@C=15pt{L \ar@{~>}[r] & M \ar@{~>}[r]
& N}$.}
\begin{demo}
Let $L,N\in\Ga$ and assume to the contrary that there exists an
infinite set of indecomposable directing modules
$\mc{M}=\{M_\lambda\}_{\lambda\in\Lambda}$ such that, for each
$\lambda\in\Lambda$, there is a path
$\xymatrix@1@C=15pt{L \ar@{~>}[r] & M_\lambda \ar@{~>}[r] & N}$
in $\ind A$. Since $\mc{M}$ is infinite and $\Ga$ has only
finitely many non-periodic $\tau$-orbits by
\cite[(2.3)]{Skow94II}, there exists an orbit $\mc{O}$ of $\Ga$
with $|\mc{O}\cap \mc{M}|=\infty$. Let $M\in \mc{O}$ and assume
without loss of generality that $\tau^mM\in \mc{M}$ for infinitely
many $m\geq 0$. Then, $M$ is left stable. Let $_l\Ga$ be the
connected component of the left stable part of $\Ga$ containing
$M$. It then easily follows from \cite[(1.4)]{CS96} that $_l\Ga$
contains no cycle and $_l\Ga$ has only finitely many
$\tau$-orbits. Then, $_l\Ga$ admits a section $\Delta$ such that
$_l\Ga$ is isomorphic to a full subquiver of $\mathbb{Z}\Delta$,
and is closed under predecessors by paths of irreducible morphisms
(see \cite[(3.4)]{L93III}). Moreover, for any predecessors $Q,Q'$
of $\Delta$, there exist at most finitely many integers $n\geq 0$
such that $Q$ is a predecessor of $\tau^nQ'$. However, since there
exists $s\geq 0$ such that $\tau^mM$ is a predecessor of $\Delta$
for all $m\geq s$, and since $\Ga$ is generalized standard and
convex, $L$ and $\tau^sM$ are two predecessors of $\Delta$ such
that $L$ is a predecessor of infinitely many $\tau^mM$, with
$m\geq s$, a contradiction. \cqfd
\end{demo}

\fait{prop}{prop dirmod}{Let $A$ be an almost laura algebra. Then
$A$
satisfies the following equivalent conditions:%
\begin{enumerate}
\item[\emph{(a)}] $\ind A \setminus (\la \cup \ra)$ contains only finitely many
directing $A$-modules.
\item[\emph{(b)}] There are only finitely many indecomposable
directing $A$-modules $M$ with a path $\xymatrix@1@C=15pt{I
\ar@{~>}[r] & M \ar@{~>}[r] & P}$ in $\ind A$ where $I$ is an
injective module and $P$ a projective module.
\item[\emph{(c)}] There are only finitely many indecomposable
directing $A$-modules $M$ with a path $\xymatrix@1@C=15pt{L
\ar@{~>}[r] & M \ar@{~>}[r] & N}$ in $\ind A$ where $L\notin \la$
and $N\notin \ra$.
\end{enumerate}}
\begin{demo}
We first show the equivalence of statement (a), (b) and (c).\\%
(a) implies (b). This follows from the fact that any injective
module (or projective module) has only finitely many successors
(or predecessors) lying in $\la$ (or in $\ra$, respectively) by
\cite[(1.5)]{AC03} (see also \cite[(3.2.6)]{LS06}).\\
(b) implies (c). This follows from (\ref{AC03 (1.6)}).\\
(c) implies (a). Assume that $\ind A \setminus (\la \cup \ra)$
contains an infinite set $(M_\lambda)_{\lambda\in\Lambda}$ of
directing modules. The set of trivial paths
$\xymatrix@1@C=15pt{M_\lambda \ar[r]^{\mbox{id}} & M_\lambda
\ar[r]^{\mbox{id}} & M_\lambda}$ contradicts (c).

\noindent Now, assume that $A$ is an almost laura algebra not
satisfying the condition (b).  Then, there exist an injective $I$,
a projective $P$ and infinitely many directing modules $M$ lying
on a path $\xymatrix@1@C=15pt{I \ar@{~>}[r] & M \ar@{~>}[r] & P}$.
By (\ref{prop almost laura}) and (\ref{lem paths}), all these
modules, including $I$ and $P$, belong to a unique component $\Ga$
of $\GA$. By (\ref{prop nsr}), $\Ga$ is generalized standard and
convex. This contradicts (\ref{lem std gen}). \cqfd
\end{demo}

As a consequence, we get the following theorem:

\fait{thm}{thm laura}{The following are equivalent for an almost
laura algebra $A$.
\begin{enumerate}
\item[\emph{(a)}] $A$ is laura.
\item[\emph{(b)}] $\ind A \setminus (\la \cup \ra)$ contains only finitely many
non-directing modules.
\item[\emph{(c)}] Any non-semiregular component of $\GA$ is
almost directed.
\item[\emph{(d)}] Any non-semiregular component of $\GA$
is quasi-directed.
\end{enumerate}}
\begin{demo}
(a) implies (b). This is obvious.\\%
(b) implies (d). Assume that $\Ga$ is a non-semiregular component
of $\GA$ and $M$ is a non-directing module in $\Ga$. By
(\ref{directed}), $M\in \ind A \setminus (\la \cup \ra)$ and the
claim follows from the assumption and (\ref{prop nsr}).\\%
(d) implies (c). This is obvious.\\%
(c) implies (a). Assume that $A$ is not laura. So $\ind A
\setminus (\la \cup \ra)$ is infinite and, by $(\ref{AC03
(1.6)})$, there exist an injective module $I$, a projective module
$P$ and infinitely many modules $M$ lying on a path
$\xymatrix@1@C=15pt{I \ar@{~>}[r] & M \ar@{~>}[r] & P}$. By
assumption, we may assume that these modules are directing. Since
$A$ is almost laura, it follows from (\ref{prop almost laura}) and
(\ref{lem paths}) that all these modules, including $I$ and $P$,
belong to the same component $\Ga$ of $\GA$. By (\ref{prop nsr}),
$\Ga$ is generalized standard and convex, which contradicts
(\ref{lem std gen}). \cqfd
\end{demo}

We get a similar characterization of almost laura algebras which
are weakly shod. Recall from \cite{CL02} that an algebra $A$ is
\textbf{weakly shod} if and only if it is laura and none of the
non-semiregular components of $\GA$ contains cycles. Moreover, a
non-semiregular component $\Ga$ is \textbf{pip-bounded} if there
exists an $n_0$ such that any path of non-isomorphisms in $\ind A$
from an injective module in $\Ga$ to a projective module in $\Ga$
has length at most $n_0$.

\fait{prop}{prop weakly shod}{The following are equivalent for an
almost laura algebra $A$.
\begin{enumerate}
\item[\emph{(a)}] A is weakly shod.
\item[\emph{(b)}] $\ind A \setminus (\la \cup \ra)$ contains only directing modules.
\item[\emph{(c)}] Any non-semiregular component of $\GA$
is directed.
\item[\emph{(d)}] Any non-semiregular component of $\GA$
is pip-bounded.
\end{enumerate}}
\begin{demo}
(a) implies (c). This follows from the above discussion. \\%
(c) implies (d). This follows from (\ref{prop nsr})
and \cite[(4.2.6)]{LS06} (see also \cite[(3.12)]{Smith04}). \\%
(d) implies (b). Assume that $M$ is a non-directing module in
$\ind A \setminus (\la \cup \ra)$.  By (\ref{AC03 (1.6)}), there
exists a path
$\xymatrix@1@C=15pt{I \ar@{~>}[r] & M \ar@{~>}[r] & P}$
in $\ind A$ for some injective module $I$ and projective module
$P$.  Since $A$ is almost laura, the modules $I,M$ and $P$ belong
to the same component $\Ga$ of $\GA$, which is therefore
non-semiregular. Obviously,
$\Ga$ is not pip-bounded, a contradiction.\\%
(b) implies (a). By (\ref{thm laura}), $A$ is laura.  Now, assume
that $\Ga$ is a non-semiregular component of $\GA$ containing a
non-directing module $M$.  By (\ref{paths to cycles}), there exist
an indecomposable injective $I$, a projective module $P$ and a
path
$\xymatrix@1@C=15pt{I \ar@{~>}[r] & M \ar@{~>}[r] & P}$.
By non-sectionality of cycles \cite{BS83,IT84} and (\ref{AC03
(1.6)}), we get $M\in \ind A \setminus (\la \cup \ra)$, a
contradiction.  So $A$ is weakly shod. \cqfd
\end{demo}

The preceding results provide new characterizations for laura and
weakly shod algebras.  We need one further lemma.

\fait{lemma}{lem LS}{Let $A$ be an algebra such that $\ind A
\setminus (\la \cup \ra)$ contains only finitely many
non-directing modules.  Then $A$ is almost laura.}
\begin{demo}
Assume that $A$ is not almost laura. Then, there exist
$L,N\notin\la\cup\ra$ such that $\rad^\infty(L,N)\neq 0$. Invoking
\cite[(4.2.2)]{LS06}, there exist infinitely many non-directing
modules $M_\lambda$ lying on a path from $L$ to $N$. Since $\la$
is closed under predecessors and $\ra$ is closed under successors,
we have $M_\lambda\notin\la\cup\ra$ for any $\lambda$. This
contradicts our assumption, and so $A$ is almost laura. \cqfd
\end{demo}

We get the following result whose proof follows from (\ref{thm
laura}), (\ref{prop weakly shod}) and (\ref{lem LS}).

\fait{cor}{cor laura and ws}{ Let $A$ be an algebra. %
\begin{enumerate}
\item[\emph{(a)}] $A$ is laura if and only if $\ind A
\setminus (\la \cup \ra)$ contains only finitely many
non-directing modules.
\item[\emph{(b)}] $A$ is weakly shod if and only if $\ind A
\setminus (\la \cup \ra)$ contains only directing modules. \cqfd
\end{enumerate}}
%

%

\subsection{Left glued algebras revisited}
    \label{Left glued algebras revisited}

A particular class of laura algebras is given by the so-called
left (or right) glued algebras.  Recall from \cite{AC94,AC03} that
an algebra $A$ is called \textbf{left glued} if the set $\ind
A\setminus \ra$ is finite.  The \textbf{right glued} algebras are
defined dually. The origin of their names comes from the fact
that, roughly speaking, the AR-quiver of any left glued algebra is
obtained by "gluing", on the left-hand-side of the AR-quiver of a
representation-finite algebra, some AR-components (without
injectives) arising from tilted algebras (see \cite{AC94} for
details).

It is well-known that left (or right) glued algebras are
characterized by the existence, in their AR-quiver, of a faithful
$\pi$-component (or $\iota$-component respectively). Recall from
\cite{C93} that an AR-component $\Ga$ is called a
\textbf{$\pi$-component} (or a \textbf{$\iota$-component})
provided all but finitely many modules in $\Ga$ are directing and
lie in the $\tau$-orbit of a projective (or an injective,
respectively).
We refer to \cite{ACLST05,Liu03} for more details concerning left
(or right) glued algebras.

The aim of this section is to show that, although laura and almost
laura algebras differ from many points of view, the "left glued"
and "right glued" versions for almost laura algebras coincide with
the usual left and right glued algebras arising from laura
algebras.

\fait{thm}{thm glued}{Let $A$ be an algebra.
\begin{enumerate}
\item[\emph{(a)}] $A$ is left glued if and only if $\rad^\infty(M,N)=0$
for all $M,N\in\ind A\setminus\ra$.
\item[\emph{(b)}] $A$ is right glued if and only if $\rad^\infty(M,N)=0$
for all $M,N\in\ind A\setminus\la$.
\end{enumerate}}
\begin{demo}
We only prove (a) since the proof of (b) is dual.\\
(a). The necessity clearly follows from the definition of left
glued algebras and \cite[(1.1)]{Smith04}, for instance.
Conversely, assume that $\rad^\infty(M,N)=0$ for all $M,N\in\ind
A\setminus\ra$.  If $\ind A=\ra$, then there is nothing to prove.
Otherwise, let $M\in\ind A\setminus\ra$ and $\Ga$ be the
AR-component containing $M$.  We show that $\Ga$ is a faithful
$\pi$-component. Let $P$ be an indecomposable projective module
such that $\Hom(P,M)\neq 0$. Since $M\notin \ra$, we have
$P\notin\ra$.  It then follows from our hypothesis that
$\rad^\infty(P,M)=0$, and so $P$ lies in $\Ga$.  So $\Ga$ contains
projective modules.  We claim that $\Ga$ contains all projective
modules.  Indeed, if this is not that case, then there exist a
projective module $P$ in $\Ga$ and a projective module $P'$ not in
$\Ga$ such that $\rad^\infty(P,P')\neq 0$ or
$\rad^\infty(P',P)\neq 0$.  Assume that $\rad^\infty(P,P')\neq 0$.
Then, since there are only finitely many predecessors of $P'$
lying in $\ra$ by \cite[(1.5)]{AC03} and \cite[(3.2.6)]{LS06}, it
follows from \cite[(1.1)]{Smith04}, for instance, that there
exists a predecessor $N$ of $P'$ such that $N\notin\ra$ but
$\rad^\infty(P,N)\neq 0$, which contradicts our hypothesis. The
same argument shows that $\rad^\infty(P',P)\neq 0$.  So $\Ga$
contains all indecomposable projective modules. In particular,
$\Ga$ is faithful.
Moreover, we have $\rad^\infty(-, \Ga)=0$. Indeed, assume that
$\rad^\infty(M', N')\neq 0$ for some indecomposables $M', N'$ with
$N'\in\Ga$.  Then, invoking (\ref{pred proj}), and recalling that
there exist only finitely many predecessors of a projective module
in $\ra$, there exists a projective module $P''$ in $\Ga$ and an
indecomposable module $M''\notin\ra$ such that $\rad^\infty(M'',
P'')\neq 0$.  This contradicts our assumption. Hence
$\rad^\infty(-,\Ga)=0$, and $\Ga$ is a $\pi$-component by
\cite[(2.1)-(2.3)]{Liu03}. Since $\Ga$ is also faithful, then $A$
is left glued. \cqfd
\end{demo}

%
%
\section{Supported almost laura algebras}
    \label{Supported almost laura algebras}

As pointed out in the discussion following (\ref{ex almost
laura}), the AR-quivers of almost laura algebras usually have many
non-semiregular components.  It is also easy to construct examples
of almost laura algebras having multicoils (in the sense of
\cite{AS92}).  With this in mind, it seems that the general shape
of the AR-quiver of an almost laura algebra is not easy to
describe.  In this section, we propose to study the AR-quiver of
left (or right) supported almost laura algebras \cite{ACT04}.

Informally, left (or right) supported algebras $A$ are those whose
left (or right) part "behaves well". For instance, any strict
laura algebra is left and right supported by \cite[(4.4)]{ACT04}.
This is however not true for almost laura algebras, as we will
see, and this additional assumption will be very useful in our
attempt to describe their AR-quivers. The main result of this
section is an analogue to the results of \cite[(3.1)]{RS04} and
\cite[(4.2.5)]{LS06}  for laura algebras and states that if $A$ is
left (or right) supported, then $A$ is almost laura if and only if
its AR-quiver has a generalized standard, convex and faithful
component (see (\ref{thm almost laura})).

Here, we recall basic features needed in the subsequent
developments.  For a full account, we refer to
\cite{ACT04,ACLST05}. By \cite{AS80}, a full subcategory $\mc{C}$
of $\mod A$ is \textbf{contravariantly finite} if for any
$N\in\mod A$, there exists a morphism
$\mor{f_{\mc{C}}}{M_{\mc{C}}}{N}$, with $M_{\mc{C}}\in\mc{C}$,
such that any morphism $\mor{f}{M}{N}$, with $M\in\mc{C}$, factors
through $f_{\mc{C}}$. The dual notion is that of a
\textbf{covariantly finite} subcategory.  Following \cite{ACT04},
an artin algebra $A$ is called \textbf{left supported} in case
$\add \la$ is contravariantly finite in $\mod A$. We define dually
the \textbf{right supported} algebras.
%
In what follows, the dual statements for right supported algebras
hold as well. We refrain from stating them.

In order to have a better description of left supported algebras,
we define, following \cite{ACT04}, two subclasses of $\la$ :
$$\begin{array}{c}
      \mc{E}_1=\{\ M\in \la \ | \ \mbox{there exists an injective } I \mbox{ and a path of irreducible } \hfill \\ %
      \qquad \mbox{ morphisms } \xymatrix@1@C=15pt{I\ar@{~>}[r]&M}\}, \mbox{ and} \hfill \\
      \mc{E}_2=\{\ M\in \la\setminus\mc{E}_1 \ | \ \mbox{there exists a projective } P\notin \la \mbox{ and  a path of } \hfill \\ %
      \qquad \mbox{ irreducible morphisms } \xymatrix@1@C=15pt{P\ar@{~>}[r]&
      \tau^{-1}M}\}.\hfill
    \end{array}$$
Moreover, we set $\mc{E}=\mc{E}_1\cup\mc{E}_2$. We also denote by
$E$ the direct sum of all indecomposable $A$-modules lying in
$\mc{E}$ and by $F$ the direct sum of a full set of
representatives of the isomorphism classes of indecomposable
projective $A$-modules not lying in $\la$. Finally, we set $T = E
\oplus F$.
The following summarizes some characterizations of left supported
algebras, as stated and proved in \cite[(Theorem A)]{ACT04} and
\cite[(Section 8)]{A07}.

\fait{thm}{thm left supported}{Let $A$ be an algebra.  The
following are equivalent:
\begin{enumerate}
\item[\emph{(a)}] $A$ is left supported.
\item[\emph{(b)}] $\add \la$ coincides with the set $\emph{Cogen}E$ of $A$-modules cogenerated by $E$.
\item[\emph{(c)}] $T = E\oplus F$ is a tilting $A$-module.
\item[\emph{(d)}] Every morphism $\mor{f}{M}{N}$ in $\ind A$, with
$M\in\la$ and $N\notin\la$ factors through $\add E$. %
\cqfd
\end{enumerate}}

\fait{remark}{Almost laura not spg}{\emph{%
Strict almost laura algebras are not left supported in general.
Indeed, for the almost laura algebra of (\ref{ex almost
laura})(b), it is easily verified that
$T= \substack{44 \\ 3}
\oplus \ 4 \ \oplus
\substack{5 \\ 24 \\ 1}$.  Since $T$ admits less indecomposable
direct summands than the number of non-isomorphic simple modules,
$T$ is not a tilting module. So $A$ is not left supported by the
above theorem.}}

We begin the study of left supported almost laura algebras with
the following lemma.  In the sequel, we write $M\fact N$ to
express that an $A$-module $M$ is a direct summand of an
$A$-module $N$.
%

\fait{lemma}{lem epsilon}{Let $A$ be an almost laura algebra. If
$M\fact T$, then the component containing $M$ also contains
injective modules.}
\begin{demo}
If $M\in \mc{E}_1$, this is clear.
If $M\in \mc{E}_2$, then there is a projective module $P\notin
\la$ and a path of irreducible morphism
$\xymatrix@1@C=15pt{P\ar@{~>}[r] & \tau^{-1}M}$.  Since $P\notin
\la$, it follows from (\ref{AC03 (1.6)}) that there is an
injective module $I$ and a path $\xymatrix@1@C=15pt{I\ar@{~>}[r] &
P}$. Since $A$ is almost laura, $I, P$ and $N$ belong to the same
component of $\GA$ by (\ref{prop almost laura}).
Finally, if $M\fact F$, then $M$ is a projective module not in
$\la$. A repetition of the above argument leads to the result.
\cqfd
\end{demo}

As a consequence, we obtain the following very useful result.

\fait{prop}{prop T}{Let $A$ be a left supported almost laura algebra.
\begin{enumerate}
\item[\emph{(a)}] If $A$ is quasitilted, then $A$ is
tilted and there exists a connecting component $\Ga$ of $\GA$
containing every indecomposable direct summand of $T$. In
particular, $\Ga$ is faithful.
\item[\emph{(b)}] If $A$ is not quasitilted, then $\GA$ has a
unique non-semiregular component $\Ga$.  Moreover, $\Ga$ contains
every indecomposable direct summand of $T$ and is faithful.
\end{enumerate}}
\begin{demo}
(a). If $A$ is quasitilted, then $A$ is tilted having $\mc{E}$ as
complete slice by \cite[(3.8)]{Smith04}. Since $F=0$ in this case,
the result follows at once.\\
(b). If $A$ is not quasitilted, let $\Ga$ be a non-semiregular
component (see (\ref{prop nsr})(b)). Then, $T$ admits an
indecomposable direct summand in $\Ga$.  Indeed, let $P$ be a
projective module in $\Ga$. If $P\notin\la$, then $P\fact F$, and
we are done. Otherwise, $P\in\la$, and since $\Ga$ contains
injective modules, we have $\Ga\cap\mc{E}\neq\emptyset$ by
\cite[(3.5)]{ACT04}. We now show that $\Ga$ contains all
indecomposable direct summands of $T$. Indeed, if this is not the
case, then there exists such a summand $T'$ of $T$ with
$\rad^\infty(\Ga,T')\neq 0$ or $\rad^\infty(T', \Ga)\neq 0$ (since
$\End_AT$ is connected). Since the component containing $T'$
contains injective modules by (\ref{lem epsilon}), we have
$\rad^\infty(T',\Ga)=0$ by (\ref{lem paths}).  So
$\rad^\infty(\Ga, T')\neq 0$. Applying (\ref{prop Hom}), we get
$T'\in\ra\setminus\la$, and so $T'\fact F$.  But then $T'$ is
projective and we get a contradiction to (\ref{lem paths}). This
proves our claim. Finally, $\Ga$ is faithful since so is $T$.
\cqfd
\end{demo}

\fait{cor}{cor capcup}{Let $A$ be a left supported almost laura
algebra. Assume that $\Ga$ is a non-semiregular
component of $\GA$ and $M\in\ind A$.%
\begin{enumerate}
\item[\emph{(a)}] $\la\cap\ra$ is finite and lies in $\Ga$.
\item[\emph{(b)}] If $M\notin \la\cup\ra$, then $M\in\Ga$.
\item[\emph{(c)}] If $M\notin\Ga$, then $M\in\la\setminus\ra$ or $M\in\ra\setminus\la$.
\end{enumerate}}
\begin{demo}
(a). Let $M\in \la\cap\ra$, and assume that $M\notin \Ga$. Since
$M\in \mbox{Cogen}E$ and $\mc{E}\subseteq \Ga$, we have
$\Hom_A(M,\Ga)\neq 0$.  By (\ref{prop Hom}), we obtain $M\notin
\ra$, a contradiction.
Now, assume to the contrary that $\la\cap\ra$ is infinite. Since
$\Ga$ has only finitely many non-periodic $\tau$-orbits by
\cite[(2.3)]{Skow94II}, there exists a $\tau$-orbit $\mc{O}$ of
$\Ga$ such that $|\mc{O}\cap(\la\cap\ra)|=\infty$. Let $M\in
\mc{O}$ and assume, without loss of generality, that $\tau^mM\in
\la\cap\ra$ for infinitely many $m\leq 0$. Then, $M$ is right
stable and it follows from \cite[(1.1)]{CL02} that there exists
$n\leq 0$ such that $\tau^nM$ is a successor of an injective
module in $\Ga$. By (\ref{AC03 (1.6)}), we have
$\tau^{n-1}M\notin\la$.  But this contradicts our assumption on
$M$. So $\la\cap\ra$ is finite.\\
(b). This follows from (\ref{prop semiregular})(a).\\
%
(c). This follows from (a) and (b). \cqfd
\end{demo}

This yields the following structure results.

\fait{lemma}{lem Hom}{Let $A$ be a left supported almost laura
algebra. Assume that $\Ga$ is a non-semiregular component of
$\GA$. Let $M\in \ind A$.  If $M\notin \Ga$, then
\begin{enumerate}
\item[\emph{(a)}] $\Hom_A(M, \Ga)\neq 0$ if and only if $M \in\la\setminus\ra$.
\item[\emph{(b)}] $\Hom_A(\Ga, M)\neq 0$ if and only if $M\in\ra\setminus\la$.
\item[\emph{(c)}] Either $\Hom_A(M, \Ga)\neq 0$ and $\Hom_A(\Ga,
M)=0$, or $\Hom_A(M, \Ga)= 0$ and $\Hom_A(\Ga, M)\neq 0$.
\end{enumerate}}
\begin{demo}
(a). Since the necessity follows from (\ref{prop Hom}), assume
that $M\in\la\setminus\ra$. Since $M\in\la\subseteq \mbox{Cogen}E$
and $\mc{E}\subseteq \Ga$, we have $\Hom_A(M, \Ga)\neq 0$. \\%
(b). Since the necessity follows from (\ref{prop Hom}), assume
that $M\in\ra\setminus\la$. Let $P$ be an indecomposable
projective module such that there exists a non-zero morphism
$\mor{\pi}{P}{M}$. If $P\in\la$, then $\pi$ factors through $\add
E$ by (\ref{thm left supported}) and so $\Hom_A(\Ga, M)\neq 0$
since $\mc{E}\subseteq \Ga$ by (\ref{prop T}). Otherwise, $P\fact
F$, and so $P\in\Ga$. Consequently, $\Hom_A(\Ga, M)\neq 0$.\\
(c). By (\ref{cor capcup}), we have $M\in\la\setminus\ra$ or
$M\in\ra\setminus\la$.  The result then follows from (a) and (b).
\cqfd
\end{demo}

\fait{thm}{thm Hom2}{Let $A$ be a left supported almost laura
algebra. Assume that $\Ga$ is a non-semiregular component of
$\GA$. Let $\Ga'\neq \Ga$ be a component of $\GA$.
\begin{enumerate}
\item[\emph{(a)}] $\Hom_A(\Ga', \Ga)\neq 0$ if and only if $\Ga' \subseteq\la\setminus\ra$.
\item[\emph{(b)}] $\Hom_A(\Ga, \Ga')\neq 0$ if and only if $\Ga'\subseteq\ra\setminus\la$.
\item[\emph{(c)}] Either $\Hom_A(\Ga', \Ga)\neq 0$ and $\Hom_A(\Ga,
\Ga')=0$, or $\Hom_A(\Ga', \Ga)= 0$ and $\Hom_A(\Ga, \Ga')\neq 0$.
\end{enumerate}
In particular, $\Ga$ is the unique faithful component of $\GA$.}
\begin{demo}
(a). Since the necessity follows from (\ref{prop Hom}), assume
that $\Ga'\subseteq\la\setminus\ra$. Let $M\in\Ga'$. By (\ref{lem
Hom}), we have $\Hom_A(M, \Ga)\neq 0$ and so $\Hom_A(\Ga',
\Ga)\neq 0$.   \\%
(b). The proof is similar to that of (a) and is left to the
reader.\\
(c). Let $M\in\Ga'$. By (\ref{lem Hom}), we have $\Hom_A(\Ga',
\Ga)\neq 0$ or $\Hom_A(\Ga, \Ga')\neq 0$.  The result then follows
from (a) and (b).\\
Finally, observe that $\Ga$ is faithful by (\ref{prop T}) and that
if $\Ga'$ was another faithful component, then we would have
$\Hom_A(\Ga,\Ga')\neq 0$ and $\Hom_A(\Ga',\Ga)\neq 0$. \cqfd
\end{demo}

\fait{remark}{trisection}{ \emph{%
Under the assumptions of (\ref{thm Hom2}) the component $\Ga$
induces a trisection in the family of AR-components (in the sense
of \cite{PS07}) : there are the components lying in
$\la\setminus\ra$, those lying in $\ra\setminus\la$ and $\Ga$.
Also, any component $\Ga'$ in $\la\setminus\ra$ maps non-trivially
to $\Ga$, which maps non-trivially to any component $\Ga''$ in
$\ra\setminus\la$. In addition, with these notations, it follows
from (\ref{thm left supported})(d) and (\ref{prop T}) that any
morphism from $\Ga'$ to $\Ga''$ factors through $\Ga$. Moreover,
by \cite[(5.5)]{ACT04}, any component lying in $\la\setminus\ra$
has no injectives and is either a postprojective component, a
semiregular tube, a component of the form $\mathbb{Z}A_{\infty}$
or a ray extension of $\mathbb{Z}A_{\infty}$. Numerous important
families of algebras accept a trisection of its module category,
notably the tilted algebras, the quasitilted algebras, the weakly
shod algebras and the laura algebras.}}

We can now prove the main result of this section, which is a
characterization of left supported almost laura algebras.

\fait{thm}{thm almost laura}{Let $A$ be a left supported algebra.
Then $A$ is almost laura if and only if $\Ga(\mod A)$ has a
generalized standard, convex and faithful component.}
\begin{demo}
The necessity follows from (\ref{prop T}), (\ref{prop nsr}) and
the fact that any connecting component is generalized standard and
convex. Conversely, assume that $\Ga$ is a generalized standard,
convex and faithful component in $\Ga(\mod A)$.  In addition,
assume that $\xymatrix@1@C=15pt{I\ar@{~>}[r] & P}$ is a path in
$\ind A$, with $I$ injective and $P$ projective. Since $\Ga$ is
faithful, there exist $M,N\in\Ga$ and a path of the form
$\xymatrix@1@C=15pt{M\ar[r] & I\ar@{~>}[r] & P \ar[r] & N}$
Since $\Ga$ is convex, then every module on this path belongs to
$\Ga$. Now, $\Ga$ being generalized standard, this path is finite.
So $A$ is almost laura by (\ref{prop almost laura}).
\cqfd
\end{demo}

At this point, we stress that the assumption of being left
supported was unnecessary to prove the sufficiency.  We then
deduce the following corollary.

\fait{cor}{cor almost laura}{Let $A$ be an algebra and assume that
$\Ga$ is a generalized standard and convex component of $\GA$. The
algebra $B=A/\ann \Ga$ is almost laura, where $\ann \Ga=\{a\in A \
| \ aM=0 \mbox{ for each }M\in\Ga\}$.}
\begin{demo}
Clearly $\Ga$ is a faithful component of $\Ga(\mod B)$.  In
addition, since $\mod B$ is a full subcategory of $\mod A$, then
$\Ga$ is generalized standard and convex as a component of
$\Ga(\mod B)$.  The result then follows from (\ref{thm almost
laura}). \cqfd
\end{demo}

The above corollary shows the importance of identifying the
generalized standard and convex components.  In the vein of
\cite{Smith04,LS06}, we then state the following result whose
proof, left to the reader, easily follows using (\ref{pred proj}).

\fait{prop}{std gen and convex}{Let $A$ be an algebra and $\Ga$ be
a component in $\GA$.  Then $\Ga$ is generalized standard and
convex if and only if any path connecting two modules in $\Ga$ is
finite.
In addition, %
\begin{enumerate}
\item[\emph{(a)}] If $\Ga$ is non-semiregular, then this is the case if
and only if any path from an injective in $\Ga$ to a projective in
$\Ga$ is finite.
\item[\emph{(b)}] If $\Ga$ is semiregular, then this is the case if and
only if any cycle $\xymatrix@1@C=15pt{M\ar@{~>}[r]&M}$, with
$M\in\Ga$, is finite. Moreover, %
\begin{enumerate}
\item[\emph{(i)}] if $\Ga$ contains injectives but no projectives, then this
occurs if and only if any path from an injective in $\Ga$ to a
module in $\Ga$ is finite;
\item[\emph{(ii)}] if $\Ga$ contains projectives but no injectives, then this
occurs if and only if any path from a module in $\Ga$ to a
projective in $\Ga$ is finite. \cqfd
\end{enumerate}
\end{enumerate}}

If $A$ is strict almost laura, then the generalized standard,
convex and faithful component of (\ref{thm almost laura}) is
non-semiregular. Since, by \cite[(3.1)]{RS04}, an algebra $A$
which is not quasitilted is laura if and only if $\GA$ has a
non-semiregular faithful and quasi-directed component, this
motivates the following problem.

\fait{prob}{prob3}{Let $A$ be a left supported strict almost laura
algebra and $\Ga$ be the unique non-semiregular component of
$\GA$.  Is $\Ga$ almost directed?}

Since strict laura algebras are left and right supported, a
positive answer would show that, for a
strict almost laura algebra $A$, the following are equivalent:%
\begin{enumerate}
\item[(a)] $A$ is left supported.
\item[(b)] $A$ is right supported.
\item[(c)] $A$ is laura.
\end{enumerate}

We end this section with a discussion of the case where $\la$ is
finite, that is contains only finitely many objects.  

\fait{prop}{prop pi}{Let $A$ be an almost laura algebra such that
$\la$ is finite. Then $\GA$ admits a faithful non-semiregular
$\pi$-component $\Ga$. In particular, $\rad^\infty(-,\Ga)=0$ and
$A$ is left glued.}
\begin{demo}
We can clearly assume that $A$ is representation-infinite.
Moreover, observe that $A$ is left supported since $\la$ is
finite, and let $\Ga$ be the (faithful) component of (\ref{prop
T}). Since $\la$ is finite and $\Ga$ is generalized standard, we
have $\rad^\infty(-,\Ga)=0$ by (\ref{thm Hom2}). In particular,
$\Ga$ contains projective modules, and so $\Ga$ is non-semiregular
by (\ref{prop semiregular}). Then $\Ga$ is a $\pi$-component by
\cite[(2.1)-(2.3)]{Liu03}. Hence $A$ is left glued. \cqfd
\end{demo}

\fait{prop}{prop repfin}{Let $A$ be an almost laura algebra. Then
$\la$ and $\ra$ are finite if and only if $A$ is
representation-finite.}
\begin{demo}
It clearly suffices to prove the necessity. If $A$ is quasitilted,
then there is nothing to show since $\ind A=\la\cup\ra$ by
\cite[(II.1.13)]{HRS96}. So, let $A$ be a strict almost laura
algebra and $\Ga$ be a non-semiregular component as in (\ref{prop
nsr})(b). By (\ref{prop pi}) and its dual, we have
$\rad^\infty(-,\Ga)=0=\rad^\infty(\Ga,-)$. So $\rad^\infty(\mod
A)=0$ and $A$ is representation-finite by \cite[(V.7.7)]{ARS97}.
\cqfd
\end{demo}

%
%

\section{Full subcategories, split-by-nilpotent extensions and
skew group algebras}
    \label{Full subcategories, split-by-nilpotent extensions and
skew group algebras}

Starting with an algebra $A$, it is frequent in the representation
theory of artin algebras to consider natural constructions giving
rise to a new algebra $B$. It is then natural to ask which
properties of $\mod A$ carry over $\mod B$ and conversely. In this
final section, we consider three different such situations and
show that almost laura algebras behave well with respect to those.

\subsection{Full subcategories}
    \label{Full subcategories}

We consider the following problem. Let $A$, $B$ be artin algebras
such that $B$ is a connected full subcategory of $A$. We choose an
idempotent $e\in A$ so that $B = eAe$. Let $P=Ae$ be the
corresponding projective $A$-module. We denote by $\pres P$ the
full subcategory of $\mod A$ formed by the \textbf{$P$-presented
modules}, that is the $A$-modules $M$
for which there exists an exact sequence, of the form %
$\xymatrix@1@C=15pt{P_1 \ar[r] & P_0 \ar[r] & M\ar[r] & 0},$
with $P_0,P_1$ in $\add P$. By \cite[(II.2.5)]{ARS97}, the functor
$\Hom_A(P,-):\xymatrix@1@C=12pt{\mod A \ar[r] & \mod B}$
induces an equivalence $\pres P \cong \mod B$, under which direct
summands of $P$ correspond to the projective $B$-modules.  In
addition, by \cite[(2.1)]{AC04}, its left inverse is
$P \otimes_B -:\xymatrix@1@C=15pt{\mod B \ar[r] & \pres P
\subseteq \mod A}$,
that is if $X$ is a $B$-module, then the $A$-module $P\otimes_B X$
is $P$-presented and $\Hom_A(P,P\otimes_B X)\cong X$,
functorially.

It is shown in \cite{AC04} that $B$ is laura (or weakly shod, or
left glued) whenever so is $A$. The following enlarges this result
to almost laura algebras.

\fait{prop}{endo}{Let $A$ be an algebra and $e$ be an idempotent
in $A$ such that $B=eAe$ is connected.  If $A$ is almost laura,
then so is $B$.}
\begin{demo}
Assume that
$f:\xymatrix@1@C=15pt{X \ar[r] & Y}$
is a morphism in $\ind B$, with $X, Y \notin \lb\cup\rb$. The
functor $P\otimes_B -$ gives a morphism
%
$P\otimes_B f : \xymatrix@1@C=15pt{P\otimes_B X \ar[r] &
P\otimes_B Y}$,
where $P\otimes_B X$ and $P\otimes_B Y$ do not lie in
$\la\cup\ra$. Indeed, if, for instance, $P\otimes_B X\in
\la\cup\ra$, then $X\cong \Hom_A(P,P\otimes_B X)\in\lb\cup\rb$ by
\cite[(2.3)]{AC04}, a contradiction. Now, since $A$ is almost
laura, we have $P\otimes_B f\notin \rad^\infty(\mod A)$, and then
$f\notin \rad^\infty(\mod B)$ since
$\Hom_A(P,-):\xymatrix@1@C=15pt{\pres P \ar[r] & \mod B}$
is an equivalence.  So $B$ is almost laura. \cqfd
\end{demo}

\fait{remark}{rem full}{\emph{%
We may ask whether an artin algebra $A$ is almost laura provided
$eAe$ is almost laura for any idempotent $e\neq 1$ of $A$. The
answer is negative, and can be easily verified on the algebra of
(\ref{counter-ex}).}}

\subsection{Split-by-nilpotent extensions}
    \label{Split-by-nilpotent extensions}

We now consider another construction. Informally, if one can
roughly think of taking full subcategories as "deleting points",
the construction we now outline can be thought of as "deleting
arrows".

Let $A$ and $B$ be artin algebras and let $Q$ be a nilpotent ideal
of $A$ (that is, $Q \subseteq \textup{rad} A$). Following
\cite{AM98}, we say that $A$ is a \textbf{split-by-nilpotent
extension of $B$ by $Q$} if there exists a split surjective
algebra morphism $\xymatrix@1@C=15pt{A \ar[r]& B}$ with kernel
$Q$.
%
For instance, if $Q^2 = 0$, then the above definition coincides
with that of the trivial extension of $B$ by $Q$. Another example
is that of one-point extension. For further examples, we refer the
reader to \cite{AZ03}.

%
We consider the change of rings functors $A \otimes_B
-:\xymatrix@1@C=15pt{\mod B \ar[r] & \mod A}$ and $B \otimes_A
-:\xymatrix@1@C=15pt{\mod A \ar[r]& \mod B}$. The image of the
functor $A \otimes_B -$ in $\mod A$ is called the category of
\textbf{induced} modules.  We have the obvious natural isomorphism
$B\otimes_A A \otimes_B - \cong 1_{\mod B}$.

In is shown in \cite{AZ03} that if $A$ is laura (or weakly shod,
or left glued), then so is $B$.  The same result holds for almost
laura algebras.

\fait{prop}{prop sbne}{Let $A$ be a split-by-nilpotent extension
of $B$ by $Q$.  If $A$ is almost laura, then so is $B$.}
\begin{demo}
Assume that
$f:\xymatrix@1@C=15pt{X \ar[r] & Y}$
is a morphism in $\ind B$, with $X, Y \notin \lb\cup\rb$. The
functor $A \otimes_B -$ gives a morphism of induced indecomposable
$A$-modules
$A \otimes_B f : \xymatrix@1@C=15pt{A \otimes_B X \ar[r] & A
\otimes_B Y}$.
Moreover, $A\otimes_B X$ and $A\otimes_B Y$ do not lie in
$\la\cup\ra$ by \cite[(2.3)]{AZ03}.
Since $A$ is almost laura, we have $A \otimes_B f\notin
\rad^\infty(\mod A)$, and then $f\notin \rad^\infty(\mod B)$ since
$B\otimes_A -$ induces an equivalence between $\mod B$ and the
induced modules in $\mod A$. Thus $B$ is almost laura. \cqfd
\end{demo}

\subsection{Skew group algebras}
    \label{Skew group algebras}

The final construction we consider is that of skew group algebras.
We are mainly motivated by the fact that skew group algebras
generally retain most features from the algebras they arise,
especially concerning homological properties. The study of the
representation theory of skew group algebras was started in
\cite{RR85,dlP83}, and more recently pursued in
\cite{FR02,ALR07,DLS07}.  We recall the relevant definitions and
refer the reader to \cite{RR85,ARS97,ALR07} for details.

Let $A$ be an artin $k$-algebra and $G$ be a group with identity
$e$. We say that $G$ \textbf{acts on} $A$ if there is a function
$\xymatrix@1@C=15pt{G\times A \ar[r] & A}$,
$\xymatrix@1@C=15pt{(\sigma,a) \ar@{|->}[r] & \sigma(a)}$,
such that:
\begin{enumerate}
\item[(a)] For each $\sigma$ in $G$, the morphism
$\sigma:\xymatrix@1@C=15pt{A \ar[r] & A}$
is an algebra automorphism;
\item[(b)] $(\sigma_1 \sigma_2)(a)=\sigma_1(\sigma_2(a))$
for all $\sigma_1,\sigma_2\in G$ and $a\in A$;
\item[(c)] $e(a)=a$ for all $a\in A$.
\end{enumerate}

Such an action induces an action of $G$ on $\mod A$ as follows :
for any $M\in \mod A$ and $\sigma\in G$, let $^\sigma M$ be the
$A$-module with the additive structure of $M$ and with the
multiplication $a\cdot m=\sigma^{-1}(a)m$, for $a\in A$ and $m\in
M$.  This allows to define an automorphism
$^\sigma (-):\xymatrix@1@C=15pt{\mod A \ar[r] & \mod A}$
for each $\sigma\in G$, where
$^\sigma f:\xymatrix@1@C=15pt{^\sigma M \ar[r] & ^\sigma N}$
is defined by
$\xymatrix@1@C=15pt{m \ar@{|->}[r] & f(m)}$
for $f\in\Hom_A(M,N)$ and $m\in M$ (see \cite[(4.1)]{ALR07}).

Suppose that $G$ acts on $A$.  The \textbf{skew group algebra}
$A[G]$ has as underlying $A$-module structure the free left
$A$-module having as basis all elements in $G$, and is endowed by
the multiplication
$(a\sigma)(b\varsigma)=a\sigma(b)\sigma\varsigma$ for all $a,b\in
A$ and $\sigma,\varsigma\in G$. Observe that $A[G]$ is generally
not connected and basic, but this will not play any role in the
sequel.

The main aim of this section is to show that if $A$ is an algebra
and $G$ is a finite group acting on $A$ and such that its order is
invertible in $A$, then $A$ is almost laura if and only if so is
$A[G]$ (see (\ref{thm skew})). It is well-known that similar
results hold for tilted, quasitilted, weakly shod and laura
algebras (see \cite[(1.2)]{ALR07}).
As we shall see, the techniques used in the proof will also result
in analogue statements for algebras having nilpotent infinite
radical and cycle-finite algebras (see (\ref{cycle-finite})).

Throughout this section, we assume that $G$ is a finite group
acting on $A$ and whose order is invertible in $A$.
Then, the natural inclusion of $A$ in $A[G]$ induces the change of
rings functors
$F:=\R\otimes_A - :\xymatrix@1@C=15pt{\mod A \ar[r] & \mod\R}$
and
$H:=\Hom_{\R}(\R, -) :\xymatrix@1@C=15pt{\mod \R \ar[r] & \mod
A}$.
We recall the following useful result from \cite[(1.1)]{RR85}.

\fait{thm}{RR (1.1)}{Let $A$ and $G$ be as above. Then
\begin{enumerate}
\item[\emph{(a)}] $(F,H)$ and $(H,F)$ are two adjoint pairs of
functors.
\item[\emph{(b)}] %
\begin{enumerate}
\item[\emph{(i)}] The unit
$\varepsilon :\xymatrix@1@C=15pt{\emph{id}_{\mod A} \ar[r] & HF}$
of the adjoint pair $(F,H)$ is a section of functors.
\item[\emph{(ii)}] The counit
$\mc{\eta} :\xymatrix@1@C=15pt{FH \ar[r] & \emph{id}_{\mod \R}}$
of the adjoint pair $(F,H)$ is a retraction of functors. \cqfd
\end{enumerate}
\end{enumerate}}

We refer to \cite[(1.1)]{RR85} for the details.  Moreover, in the
sequel, we shall use the following notations. We denote by
$\phi :\xymatrix@1@C=15pt{\Hom_{\R}(\F{-}, ?) \ar[r] &
\Hom_{A}(-,\Ho{?})}$
the natural equivalence associated to the adjoint pair $(F,H)$.
%
On the other hand, we denote by
$\psi :\xymatrix@1@C=15pt{\Hom_{A}(\Ho{?}, -) \ar[r] &
\Hom_{\R}(?,\F{-})}$
the natural equivalence associated to the adjoint pair $(H,F)$.
Finally, we let $\mu$ and $\rho$ be the unit and counit of this
adjoint pair.

With these notations, we have (see \cite[(p.~118)]{M65}, for
instance) the following useful lemma.

\fait{lemma}{lem unit}{Let $M$ be an $A$-module and
$X$ be an $\R$-module.%
\begin{enumerate}
\item[\emph{(a)}] If $f\in \Hom_{\R}(F(M),X)$, then $\phi(f)=H(f)\circ
\varepsilon_M$.
\item[\emph{(b)}] If $f\in \Hom_A(M, H(X))$, then $\phi^{-1}(f) =
\eta_X\circ F(f)$.
\item[\emph{(c)}] If $f\in \Hom_A(H(X), M)$, then $\psi(f)=F(f)\circ
\mu_X$.
\item[\emph{(d)}] If $f\in\Hom_{\R}(X,F(M))$, then
$\psi^{-1}(f)=\rho_M\circ H(f)$. \cqfd
\end{enumerate}}

We recall that given two categories $\mc{C}$ and $\mc{D}$, a
functor
$\mc{F} :\xymatrix@1@C=15pt{\mc{C} \ar[r] & \mc{D}}$
is called a \textbf{radical functor} if, for any objects $M,N$ in
$\mc{C}$, we have $\mc{F}(\rad_{\mc{C}}(M,N))\subseteq
\rad_{\mc{D}}(\mc{F}(M),\mc{F}(N))$. For instance, any full
functor is radical.

\fait{prop}{prop radical}{The functors
$F$ 
and
$H$ 
are radical functors.}
\begin{demo}
We first show that $F$ is a radical functor. Let $M,N$ be
indecomposable $A$-modules and let $f\in \rad_A(M,N)$.  Now,
assume to the contrary that $F(f)\notin \rad_{\R}(F(M),F(N))$. So,
there exist an indecomposable $\R$-module $X$ together with a
section
$\iota : \xymatrix@1@C=15pt{X \ar[r] & F(M)}$
and a retraction
$\pi : \xymatrix@1@C=15pt{F(N) \ar[r] & X}$
such that the composition $\pi\circ F(f) \circ\iota$ is an
isomorphism. Denote by
$\omega$
the left inverse of $\iota$. Applying $H$ gives a commutative
diagram:
$$\xymatrix@R=20pt@C=40pt{H(X)\ar@<.5ex>[r]^{H(\iota)} & H(F(M))
\ar@<.5ex>[l]^{H(\omega)} \ar[r]^{H(F(f))} & H(F(N))
\ar[r]^{H(\pi)} & H(X)\\
& M\ar[ul]^{\phi(\omega)} \ar[u]_{\varepsilon_M} \ar[r]^f & N
\ar[u]_{\varepsilon_N} \ar[ur]_{\phi(\pi)}}$$
where the first row is an isomorphism,
$H(\omega)\circ\varepsilon_M=\phi(\omega)$ and
$H(\pi)\circ\varepsilon_N=\phi(\pi)$ by (\ref{lem unit})(a) and
$\varepsilon_N\circ f = H(F(f))\circ \varepsilon_N$ by (\ref{RR
(1.1)})(b). Since $\phi$ is a bijection and $\omega\neq 0$, we
have $\phi(\omega)\neq 0$. It then follows from the
indecomposability of $M$ that $H(\iota)\circ
\phi(\omega)=\varepsilon_M$ and so $\phi(\omega)$ is a section.
Therefore, we have
$$\phi(\pi)\circ f= H(\pi)\circ \varepsilon_N\circ f = H(\pi)\circ
H(F(f))\circ \varepsilon_M = H(\pi)\circ H(F(f))\circ
H(\iota)\circ \phi(\omega).$$
Since $ H(\pi)\circ H(F(f))\circ H(\iota)$ is an isomorphism and
$\phi(\omega)$ is a section, then $\phi(\pi)\circ f$ is a section.
In particular, $f$ is a section, a contradiction since $N$ is
indecomposable.  So $F(f)\in \rad_{\R}(F(M),F(N))$ and $F$ is a
radical functor.
Using (\ref{lem unit})(b) and the fact that $\eta$ is a retraction
of functors, one can show in a similar way that $H$ is also a
radical functor.
\cqfd
\end{demo}

Since almost laura algebras are defined in terms of the behavior
of their infinite radicals, the knowledge of each power of the
radical is rather important. As a consequence of the above
proposition, we now show that the maps $\phi$ and $\psi$ can be
used to relate the different powers of the radicals of $\mod A$
and $\mod A[G]$.

\fait{prop}{prop rad}{Let $A$ and $G$ be as above. Let $n\geq 1$,
$M$ be an $A$-module and $X$ be an $\R$-module. Then,
\begin{enumerate}
\item[\emph{(a)}] $\phi(\rad^n_{\R}(F(M),X))=\rad^n_A(M,H(X))$;
\item[\emph{(b)}] $\psi(\rad^n_{A}(H(X),M))=\rad^n_{\R}(X,F(M))$.
\end{enumerate}
}
\begin{demo}
We only prove (a) since the proof of (b) is similar.\\
(a). Assume that $f\in\rad_{\R}^n(F(M),X)$, and let $F(M)=Y_0,
Y_1, \dots, Y_n=X$ and $f_i\in\rad_{\R}(Y_{i-1}, Y_i)$, with
$1\leq i\leq n-1$, be such that $f=f_n f_{n-1}\cdots f_1$. Then,
by (\ref{lem unit})(a), we have
$$\phi(f)=H(f)\circ \varepsilon_M=H(f_n)\circ \cdots \circ
H(f_1)\circ \varepsilon_M.$$
Since $H$ is a radical functor by (\ref{prop radical}), we have
$H(f_i)\in \rad_A(H(Y_{i-1}), H(Y_i))$ for each $i$.  So
$\phi(f)\in \rad^n_A(M,H(X))$.
%
Similarly, if $h\in \rad_A^n(M, H(X))$, then
$\phi^{-1}(h)\in\rad^n_{\R}(F(M), X)$.  The result follows. \cqfd
\end{demo}

The following two corollaries are generalizations of
\cite[(4.4)]{ALR07} and \cite[(4.6)]{ALR07} respectively. But
first, we need to recall from \cite[(4.3)]{ALR07} that if $X$ is
an indecomposable $\R$-module, then there exists an indecomposable
$A$-module $M$ such that $M\fact H(X)$ and $X\fact F(M)$.

\fait{cor}{ALR (4.4)}{Let $n\geq 1$ and $M,N$ be indecomposable
$A$-modules such that $\rad_A^n(M,N)\neq 0$.%
\begin{enumerate}
\item[\emph{(a)}] For any direct summand $X$ of
$F(M)$, we have $\rad_{\R}^n(X, F(N))\neq 0$;
\item[\emph{(b)}] For any direct summand $Y$ of
$F(N)$, we have $\rad_{\R}^n(F(M), Y)\neq 0$.
\end{enumerate}}
\begin{demo}
We only prove (a) since the proof of (b) is similar.\\
(a). By \cite[(1.8)]{RR85}, we have an indecomposable
decomposition
$F(M)\cong\oplus_{i=1}^m X_i$
in $\mod\R$ such that $H(X_i)\cong\oplus_{\sigma\in G_i} \ ^\sigma
M$ for some $G_i\subseteq G$.  In addition, for each $i$, and each
$\gamma\in G$, there exists $\sigma\in G_i$ with $^\gamma M \cong
\ ^\sigma M$. In particular, we can assume that $M\fact H(X_i)$
for each $i$.
We need to show that $\rad_A^n(X_i, F(N))\neq 0$ for each $i$ and,
by (\ref{prop rad})(b), it is sufficient to show that
$\rad_A^n(H(X_i),N)\neq 0$.  Since $M$ is a direct summand of
$H(X_i)$ for each $i$, this is clearly the case. \cqfd
\end{demo}

\fait{cor}{ALR (4.6)}{Let $n\geq 1$ and $X,Y$ be indecomposable
$\R$-modules such that $\rad_{\R}^n(X,Y)\neq 0$.  Then, for all
indecomposable $A$-modules $M,N$ such that $X\fact F(M)$ and
$Y\fact F(N)$, there exists $\sigma\in G$ such that $\rad_A^n(M, \
^\sigma N)\neq 0$.}
\begin{demo}
Let $M$ and $N$ be as in the statement. Then, by hypothesis, we
have $\rad^n_{\R}(F(M),F(N))\neq 0$, and thus
$\rad^n_A(M,H(F(N)))\neq 0$ by (\ref{prop rad}).  Since on the
other hand we have
$H(F(N)) \cong \oplus_{\sigma\in G} \ ^\sigma N$
by \cite[(1.8)]{RR85}, there exists $\sigma \in G$ with
$\rad^n_A(M, \ ^\sigma N )\neq 0$. \cqfd
\end{demo}

We also get the following corollary, which complements
\cite[(4.5)(4.7)]{ALR07}.

\fait{cor}{cor paths}{%
\begin{enumerate}
\item[\emph{(a)}] Let
$\xymatrix@1@C=15pt{M_0 \ar[r]^{f_1} & M_1 \ar[r]^{f_2} & \cdots
\ar[r]^{f_t} & M_t}$
be a path in $\ind A$, with $f_i\in \rad^{n_i}_A(M_{i-1}, M_i)$
for each $i$. For any indecomposable $X_0\fact F(M_0)$, there
exists a path
$\xymatrix@1@C=15pt{X_0 \ar[r]^{g_1} & X_1 \ar[r]^{g_2} & \cdots
\ar[r]^{g_t} & X_t}$
in $\ind \R$ with $X_i \fact F(M_i)$, $M_i\fact H(X_i)$ and
$g_i\in \rad^{n_i}_{\R}(X_{i-1}, X_i)$ for each $i$.
\item[\emph{(b)}] Let
$\xymatrix@1@C=15pt{X_0 \ar[r]^{g_1} & X_1 \ar[r]^{g_2} & \cdots
\ar[r]^{g_t} & X_t}$
be a path in $\ind \R$, with $g_i\in \rad^{n_i}_{\R}(X_{i-1},
X_i)$ for each $i$. For any indecomposable $M_0$ such that $X_0
\fact F(M_0)$, there exist $\sigma_1, \sigma_2, \dots, \sigma_t\in
G$ and a path
$\xymatrix@1@C=15pt{M_0 \ar[r]^{f_1} & \ ^{\sigma_1} M_1
\ar[r]^{f_2} & \cdots \ar[r]^{f_t} & \ ^{\sigma_t} M_t}$
in $\ind A$ with $M_i \fact H(X_i)$, $X_i\fact F(M_i)$ and $f_i\in
\rad^{n_i}_{A}( \ ^{\sigma_{i-1}} M_{i-1}, \ ^{\sigma_{i}} M_i)$
for each $i$.
\end{enumerate}}
\begin{demo}
(a). Since $X_0 \fact F(M_0)$ and $\rad^{n_1}_A(M_0,M_1)\neq 0$,
it follows from (\ref{ALR (4.4)}) that
$\rad^{n_1}_A(X_0,F(M_1))\neq 0$. Hence there exists an
indecomposable $X_1\fact F(M_1)$ with $\rad^{n_1}_A(X_0,X_1)\neq
0$.  The result follows from an obvious induction.  Observe that
$M_i\fact H(X_i)$ for each $i$ by the proof of (\ref{ALR
(4.4)}).\\
(b). Let $M_0, M_1$ be indecomposable $A$-modules such that
$X_i\fact F(M_i)$, for $i=1,2$. By (\ref{ALR (4.6)}), there exists
$\sigma_1\in G$ such that $\rad_{A}^{n_1}(M_0, \
^{\sigma_1}M_1)\neq 0$. Similarly, there exists an indecomposable
$M_2$ such that $X_2\fact F(M_2)$ together with an element
$\sigma'_2\in G$ such that $\rad_A^{n_2}(M_1, \
^{\sigma'_2}M_2)\neq 0$. Applying the automorphism
$^{\sigma_1} (-):\xymatrix@1@C=15pt{\mod A \ar[r] & \mod A}$
we obtain $\rad_A^{n_2}(^{\sigma_1}M_1, \ ^{\sigma_2}M_2)\neq 0$,
where $\sigma_2=\sigma_1\sigma'_2$. The result now follows from an
obvious induction. Observe that $M_i\fact H(X_i)$ for each $i$ by
the proof of (\ref{ALR (4.4)}). \cqfd
\end{demo}

We are now ready to prove the main result of this section.

\fait{thm}{thm skew}{Let $A$ be an algebra and $G$ be a
finite group acting on $A$ and whose order is invertible in $A$. %
\begin{enumerate}
\item[\emph{(a)}] $A$ is almost laura if and only if so is $\R$.
\item[\emph{(b)}] $A$ is strict almost laura if and only if so is $\R$.
\end{enumerate}}
\begin{demo}
(a). Assume that $A$ is almost laura and let
$\xymatrix@1@C=15pt{X \ar[r]^-{g} & Y}$
be a morphism in $\ind \R$, with $X, Y \notin \lr\cup\rr$ and
$g\in\rad^{n}_{\R}(X, Y)$. By (\ref{cor paths})(b), there exist
$\sigma\in G$ and a morphism
$\xymatrix@1@C=15pt{M \ar[r]^{f} & ^\sigma N}$
in $\ind A$ with $f\in\rad^{n}_{A}(M, ^{\sigma}N)$.
%
%
In addition, by \cite[(5.1)(5.3)]{ALR07}, we have $M, ^\sigma N
\notin \la\cup\ra$.   Since $A$ is almost laura, $f$ does not
belong to $\rad^\infty(\mod A)$, and so $g$ does not belong to
$\rad^\infty(\mod \R)$.  Hence $\R$ is almost laura.
The converse is proven in the same way, using (\ref{cor paths})(a)
instead of (\ref{cor paths})(b).\\
(b). This follows from (a) and \cite[(III.1.6)]{HRS96}. \cqfd
\end{demo}

%

Our work on the infinite radical carries consequences on other
classes of algebras, for instance on cycle-finite algebras and
algebras having nilpotent infinite radical. Recall from
\cite{AS90} that an algebra $A$ is \textbf{cycle-finite} if no
cycle in $\ind A$ contains morphisms in $\rad^\infty(\mod A)$.
Examples of cycle-finite algebras are all representation-finite
algebras, tame tilted algebras \cite{R84}, tubular algebras
\cite{R84}, iterated tubular algebras \cite{PT90}, and multicoil
algebras \cite{AS92}. It is known (see \cite{AS90}) that every
cycle-finite algebra is of tame representation type.

On the other hand, given an algebra $A$, it is important to study
the nilpotency of the infinite radical of $\mod A$ in order to
understand the complexity of $\mod A$. This has been considered,
for instance, in \cite{CMMSk94,KS91,AC03}. More precisely, we say
that $\rad^\infty(\mod A)$ is \textbf{nilpotent} if there exists
an integer $n\geq 1$ such that $(\rad^\infty(\mod A))^n=0$. Such a
minimal integer $n$ is then called the \textbf{index of
nilpotency} of $\rad^\infty(\mod A)$.

%

We have the following result.

\fait{prop}{cycle-finite}{Let $A$ be an algebra and $G$ be a
finite group acting on $A$ and whose order is invertible in $A$.
\begin{enumerate}
\item[\emph{(a)}] The infinite
radical of $\mod A$ is nilpotent if and only if so is the infinite
radical of $\mod A[G]$ and, in this case, they have the same index
of nilpotency.
\item[\emph{(b)}] $A$ is cycle-finite if and only if so is $\R$. \\ %
Moreover, in this case, $A$ is domestic if and only if so is $\R$.
\end{enumerate}}
\begin{demo}
(a). Assume that there exists an integer $n\geq 1$ such that
$(\rad^\infty(\mod A))^n=0$ but $(\rad^\infty(\mod \R))^n\neq 0$
Thus, there exists a path
$\xymatrix@1@C=15pt{X_0 \ar[r]^{g_1} & X_1 \ar[r]^{g_2} & \cdots
\ar[r]^{g_n} & X_n}$
in $\ind \R$ such that $g_i\in\rad_{\R}^\infty(\mod \R)$ for each
$i$ and $g=g_n\cdots g_2g_1\neq 0$.  Now, since $H$ is faithful by
(\ref{RR (1.1)})(b)(ii) and a radical functor by (\ref{prop
radical}), we have $0\neq H(g)\in (\rad_A^\infty(\mod A))^n$, a
contradiction. So $(\rad_{\R}^\infty(\mod \R))^n=0$.
The converse is proven in the same way, using $F$ and invoking
(\ref{RR (1.1)})(b)(i) instead of (\ref{RR (1.1)})(b)(ii).
\\%
(b). Assume that $A$ is cycle-finite and let
$\xymatrix@1@C=15pt{X=X_0 \ar[r]^-{g_1} & X_1 \ar[r]^{g_2} &
\cdots \ar[r]^-{g_t} & X_t=X}$
be a cycle in $\ind \R$, with $g_i\in\rad^{n_i}_{\R}(X_{i-1},
X_i)$ for each $i$. By (\ref{cor paths})(b), there exist
$\sigma_1, \sigma_2, \dots, \sigma_t\in G$ and a path of the form
$\delta :  \xymatrix@1@C=15pt{M_0 \ar[r]^-{f_1} & \ ^{\sigma_1}M_1
\ar[r]^{f_2} & \cdots \ar[r]^{f_t} & \ ^{\sigma_t}M_t},$
in $\ind A$, with $f_i\in\rad^{n_i}_{\R}(^{\sigma_{i-1}}M_{i-1},
^{\sigma_{i}}M_i)$ for each $i$. Moreover, by \cite[(1.8)]{RR85},
we have $^{\sigma} M_0\cong M_t$ for some $\sigma\in G$ and thus
$^{\sigma\sigma_t} M_0\cong \ ^{\sigma_t}M_t$. Let
$\tau=\sigma\sigma_t$ and $m$ be the order of $\tau$ in $G$.
Applying repeatedly the functor $^\tau (-):\xymatrix@1@C=15pt{\mod
A \ar[r] & \mod A}$ on $\delta$ yields a cycle
$$\xymatrix@1@C=25pt{M_0 \ar@{~>}[r]^{\delta} & \ ^{\tau}M_0
\ar@{~>}[r]^{^{\tau}\delta} & \ ^{\tau^2}M_0
\ar@{~>}[r]^{^{\tau^2}\delta} & \cdots
\ar@{~>}[r]^-{^{\tau^m}\delta} &  \ ^{\tau^m}M_0=M_0}$$
Since $A$ is cycle-finite, no morphism in $\delta$ belongs to
$\rad^\infty(\mod A)$, and so no $g_i$ belongs to
$\rad^\infty(\mod \R)$.  Hence $\R$ is cycle-finite. \\
On the other hand, assume that $\R$ is cycle-finite and let
$$\xymatrix@1@C=15pt{M=M_0 \ar[r]^-{g_1} & M_1 \ar[r]^{g_2} &
\cdots \ar[r]^-{g_t} & M_t=M}$$
be a cycle in $\ind A$, with $g_i\in\rad^{n_i}_A(M_{i-1}, M_i)$
for each $i$. Let $F(M)=\oplus_{j=1}^m X_j$ be an indecomposable
decomposition in $\mod \R$.  Then, for each $j$, there exists by
(\ref{cor paths})(a) a path in $\ind \R$ of the form
$\delta_j : \xymatrix@1@C=15pt{X_j \ar@{~>}[r] & X_{s_j}}$
with $1\leq s_j\leq m$ containing at least one morphism in
$\rad^{n_i}_{\R}(\mod \R)$ for each $1\leq i \leq t$. Let
$\mor{s}{\{1,2,\dots,m\}}{\{1,2,\dots,m\}}$ be the application
defined by $s(j)=s_j$.  Then,
there exist $j$ and $q$ such that $j=s^q(j)$.  Consequently, there
is a cycle 
$$\xymatrix@1@C=25pt{X_j \ar@{~>}[r]^{\delta_j} & X_{s(j)}
\ar@{~>}[r]^{\delta_{s(j)}} & X_{s^2(j)}
\ar@{~>}[r]^{\delta_{s^2(j)}} & \cdots
\ar@{~>}[r]^-{\delta_{s^q(j)}} & X_{s^q(j)}=X_j}$$
containing morphisms in $\rad^{n_i}_{\R}(\mod \R)$ for each $1\leq
i \leq t$. Since $\R$ is cycle-finite, no morphism in this path
belongs to the infinite radical, and so $A$ is cycle-finite.
The latter part directly follows from (a) and \cite[(5.1)]{Sk95}.
\cqfd

\end{demo}

\fait{remark}{T-nilpotent}{\emph{%
Recall from \cite{KS91} that $\rad^\infty(\mod A)$ is called
\textbf{left (or right) T-nilpotent} if for each sequence
$(f_i)_{i\in\mathbb{N}}$ in $\rad^\infty(\mod A)$, there exists a
natural number $m$ such that $f_m \cdots f_1=0$ (or $f_1\cdots
f_m=0$, respectively). It is easily seen that the proof of
(\ref{cycle-finite})(a) can be adapted to show that
$\rad^\infty(\mod A)$ is left (or right) $T$-nilpotent if and only
if so is $\rad^\infty(\mod\R)$.}}

%
%
\bigskip
\noindent {\bf ACKNOWLEDGEMENTS.} Most of the results in this
paper are part of the author's Ph.D. thesis, completed at the
University of Sherbrooke (Canada). The author thanks his PhD
advisor Ibrahim Assem for many helpful comments, as well as
Fl{\'a}vio Coelho for his assistance regarding Theorem \ref{thm
glued}.

%
%


\begin{thebibliography}{10}

\bibitem{A07}
I.~Assem.
\newblock {L}eft sections and the left part of an artin algebra.
\newblock Preprint.

\bibitem{AC94}
I.~Assem and F.U. Coelho.
\newblock {G}lueings of tilted algebras.
\newblock {\em J. Pure Appl. Algebra}, 96(3):225--243, 1994.

\bibitem{AC03}
I.~Assem and F.U. Coelho.
\newblock {T}wo-sided gluings of tilted algebras.
\newblock {\em J. Algebra}, 269(2):456--479, 2003.

\bibitem{AC04}
I.~Assem and F.U. Coelho.
\newblock {E}ndomorphism algebras of projective modules over laura algebras.
\newblock {\em {J}. {A}lgebra {A}ppl.}, 3(1):49--60, 2004.

\bibitem{ACLST05}
I.~Assem, F.U. Coelho, M.~Lanzilotta, D.~Smith, and S.~Trepode.
\newblock {A}lgebras determined by their left and right parts.
\newblock {\em Algebraic structures and their representations}, 13--47,
{C}ontemp. {M}ath., 376, {A}mer. {M}ath. {S}oc., {P}rovidence,
{RI}, 2005.

\bibitem{ACT04}
I.~Assem, F.U. Coelho, and S.~Trepode.
\newblock {T}he left and the right parts of a module category.
\newblock {\em J. Algebra}, 281(2):518--534, 2004.

\bibitem{ALR07}
I.~Assem, M.~Lanzilotta, and M.J. Redondo.
\newblock {L}aura skew group algebras.
\newblock {\em Comm. Algebra}, 35(7):2241--2257, 2007.

\bibitem{AM98}
I.~Assem and N.~Marmaridis.
\newblock {T}ilting modules over split-by-nilpotent extensions.
\newblock {\em Comm. Algebra}, 26(5):1547--1555, 1998.

\bibitem{AS90}
I.~Assem and A.~Skowro{\'n}ski.
\newblock {M}inimal representation-infinite coil algebras.
\newblock {\em {M}anuscripta {M}ath.}, 67(3):305--331, 1990.

\bibitem{AS92}
I.~Assem and A.~Skowro{\'n}ski.
\newblock Indecomposable modules over multicoil algebras.
\newblock {\em {M}ath. {S}cand.}, 71(1):31--61, 1992.

\bibitem{AZ03}
I.~Assem and D.~Zacharia.
\newblock {O}n split-by-nilpotent extensions.
\newblock {\em {C}olloq. {M}ath.}, 98(2):259--275, 2003.

\bibitem{ARS97}
M.~Auslander, I.~Reiten, and S.O. Smal{\o}.
\newblock {\em {R}epresentation {T}heory of {A}rtin {A}lgebras}.
\newblock Cambridge Studies in Advanced Mathematics, 36. Cambridge
  University Press, Cambridge, 1997. xiv+425 pp.

\bibitem{AS80}
M.~Auslander and S.O. Smal{\o}.
\newblock {P}reprojective modules over artin algebras.
\newblock {\em J. Algebra}, 66(1):61--122, 1980.

\bibitem{BS83}
R.~Bautista and S.O. Smal{\o}.
\newblock {N}on-existent cycles.
\newblock {\em {C}omm. {A}lgebra}, 11(6):1755--1767, 1983.

\bibitem{C93}
F.U. Coelho.
\newblock {C}omponents of {A}uslander-{R}eiten quivers with only preprojective
  modules.
\newblock {\em J. Algebra}, 157(2):472--488, 1993.

\bibitem{CL02}
F.U. Coelho and M.~Lanzilotta.
\newblock {O}n non-semiregular components containing paths from injective to
  projective modules.
\newblock {\em Comm. Algebra}, 30(10):4837--4849, 2002.

\bibitem{CMMSk94}
F.U. Coelho, E.N. Marcos, H.A. Merklen, and A.~Skowro{\'n}ski.
\newblock {M}odule categories with infinite radical square zero are of finite
  type.
\newblock {\em Comm. Algebra}, 22(11):4511--4517, 1994.

\bibitem{CS96}
F.U. Coelho and A.~Skowro{\'n}ski.
\newblock {O}n {A}uslander-{R}eiten components for quasitilted algebras.
\newblock {\em Fund. Math.}, 149(1):67--82, 1996.

\bibitem{DLS07}
J.~Dionne, M.~Lanzilotta, and D.~Smith.
\newblock {S}kew group algebras of piecewise hereditary algebras are piecewise hereditary.
\newblock Preprint.

\bibitem{DS06}
J.~Dionne and D.~Smith.
\newblock {A}rticulations of algebras and their homological properties.
\newblock {\em {J}. {A}lgebra {A}ppl.}, 5(3):1--15, 2006.

\bibitem{FR02}
O.~Funes and M.~J. Redondo.
\newblock {S}kew group algebras of simply connected algebras.
\newblock {\em {A}nn. {S}ci. {M}ath. {Q}u\'ebec}, 26(2):171--180, 2002.

\bibitem{HRS96}
D.~Happel, I.~Reiten, and S.O. Smal{\o}.
\newblock {T}ilting in abelian categories and quasi-tilted algebras.
\newblock {\em Mem. Am. Math. Soc.}, 120(575):viii+ 88 pp., 1996.

\bibitem{IT84}
K.~Igusa and G.~Todorov.
\newblock {A} characterization of finite {A}uslander-{R}eiten quivers.
\newblock {\em J. Algebra}, 89(1):148--177, 1984.

\bibitem{KS91}
O.~Kerner and A.~Skowro{\'n}ski.
\newblock {O}n module categories with nilpotent infinite radical.
\newblock {\em {C}ompositio {M}ath.}, 77(3):313--333, 1991.

\bibitem{LS06}
M.~Lanzilotta and D.~Smith.
\newblock Laura algebras and quasi-directed components.
\newblock {\em {C}olloq. {M}ath.}, 105(2):179--196, 2006.

\bibitem{L93III}
S.~Liu.
\newblock {S}emi-stable components of an {A}uslander-{R}eiten quiver.
\newblock {\em {J}. {L}ondon {M}ath. {S}oc. (2)}, 47(3):405--416, 1993.

\bibitem{Liu03}
S.~Liu.
\newblock {P}reprojective modules and {A}uslander-{R}eiten components.
\newblock {\em {C}omm. {A}lgebra}, 31(12):6051--6061, 2003.

\bibitem{M65}
B.~Mitchell.
\newblock {\em {T}heory of categories}, volume XVII of {\em {P}ure and
  {A}pplied in {M}athematics}.
\newblock Academic Press, {N}ew {Y}ork - {L}ondon, 1965.

\bibitem{PT90}
J.~A. de~la Pe\~na and Tomé.
\newblock Iterated tubular algebras.
\newblock {\em {J}. {P}ure {A}ppl. {A}lgebra}, 64(3):303--314, 1990.

\bibitem{dlP83}
J.A. de~la Pe\~na.
\newblock {A}utomorfismos, \'algebras torcidas y cubiertas.
\newblock Ph.D. thesis, Universidad Nacional Auton{\'o}ma de México, 1983.

\bibitem{PS07}
J.A. de~la Pe\~na and I.~Reiten.
\newblock Trisections of module categories.
\newblock {\em Colloq. Math.}, 107(2):191--219, 2007.

\bibitem{RR85}
I.~Reiten and C.~Riedtmann.
\newblock {S}kew group algebras in the representation theory of artin algebras.
\newblock {\em J. Algebra}, 92(1):224--282, 1985.

\bibitem{RS04}
I.~Reiten and A.~Skowro{\'n}ski.
\newblock {G}eneralized double tilted algebras.
\newblock {\em {J}. {M}ath. {S}oc. {J}apan}, 56(1):269--288, 2004.

\bibitem{R84}
C.M. Ringel.
\newblock {\em {T}ame algebras and integral quadratic forms}.
\newblock Number 1099 in Lecture Notes in Mathematics. Springer-Verlag, Berlin
  Heidelberg New York Tokyo, 1984.

\bibitem{Skow94II}
A.~Skowro{\'n}ski.
\newblock {G}eneralized standard {A}uslander-{R}eiten components.
\newblock {\em {J}. {M}ath. {S}oc. {J}apan}, 46(3):517--543, 1994.

\bibitem{Sk94II}
A.~Skowro{\'n}ski.
\newblock {M}inimal representation-infinite artin algebras.
\newblock {\em {M}ath. {P}roc. {C}am. {P}hi. {S}oc.}, 116(2):229--243, 1994.

\bibitem{Sk94III}
A.~Skowro{\'n}ski.
\newblock {R}egular {A}uslander-{R}eiten components containing directing
  modules.
\newblock {\em {P}roc. {A}mer. {M}ath. {S}oc.}, 120(1):19--26, 1994.

\bibitem{Sk95}
A.~Skowro{\'n}ski.
\newblock {C}ycle-finite algebras.
\newblock {\em {J}. {P}ure {A}ppl. {A}lgebra}, 103(1):105--116, 1995.

\bibitem{Sk03}
A.~Skowro{\'n}ski.
\newblock {O}n artin algebras with almost all indecomposable modules of
  projective or injective dimension at most one.
\newblock {\em Cent. Eur. J. Math.}, 1(1):108--122, 2003.

\bibitem{Smith04}
D.~Smith.
\newblock {O}n generalized standard {A}uslander-{R}eiten components having only
  finitely many non-directing modules.
\newblock {\em J. Algebra}, 279(2):493--513, 2004.

\end{thebibliography}

\medskip

\noindent D.~Smith; Department of Mathematical Sciences, Norwegian
University of Science and Technology, N-7491 Trondheim,
Norway\\
E-mail address: david.smith@math.ntnu.no
\end{document}